\documentclass[10pt]{article}
\usepackage{amsmath,amssymb}

\newcommand{\bra}{\langle}
\newcommand{\ket}{\rangle}
\newcommand{\ep}{\hfill {$\Box$}}

\newtheorem{thm}{Theorem}[section]
\newtheorem{cor}{Corollary}[section]
\newtheorem{defin}{Definition}[section]
\newtheorem{lem}{Lemma}[section]
\newtheorem{prop}{Proposition}[section]

\title{Commutation Relations for Unitary Operators III}
\author{\small{M.A. Astaburuaga, O. Bourget, V.H. Cort\'es \footnote{Supported by the Grants Fondecyt 1080455, 1080675, 1120786, ICM PROY-P07-027-F, ECOS-Conicyt C10E10}}\\
\small{Facultad de Matem\'aticas, Pontificia Universidad Cat\'olica de Chile,}\\
\small{Av. Vicu\~na Mackenna 4860, Macul, Santiago, Chile}\\
\small{E-mail: bourget@mat.puc.cl}\\
\small{phone: (56 2) 354 4509}\\
}
\date{}

\begin{document}
\maketitle

\begin{abstract}
Let $U$ be a unitary operator defined on some infinite-dimensional complex Hilbert space ${\cal H}$. Under some suitable regularity assumptions, it is known that a local positive commutation relation between $U$ and an auxiliary self-adjoint operator $A$ defined on ${\cal H}$ allows to prove  that the spectrum of $U$ has no singular continuous spectrum and a finite point spectrum, at least locally. We prove that under stronger regularity hypotheses, the local regularity properties of the spectral measure of $U$ are improved, leading to a better control of the decay of the correlation functions. As shown in the applications, these results may be applied to the study of periodic time-dependent quantum systems, classical dynamical systems and spectral problems related to the theory of orthogonal polynomials on the unit circle.
\end{abstract}


\section{Introduction}

The spectral analysis of unitary operators defined on Hilbert spaces is a natural tool in the study of the long-time behavior of periodic time-dependent quantum systems \cite{ev}. It also appears in the theory of orthogonal polynomials \cite{S1}, \cite{S2} and the study of classical dynamical systems e.g \cite{N1}, \cite{N2}.

The commutation relations satisfied by an operator may be relevant to determine its spectral properties. This approach has been developped to a large extent for self-adjoint operators to analyze either its discrete spectrum or its essential component by means of some positive commutator methods. The development of these methods within the spectral theory of unitary operators has been historically delayed, although this gap has been now partly filled regarding the development of the positive commutator theory \cite{abcf}, \cite{frt} and \cite{abc}.

This manuscript is focused on the relationships between the existence of a positive commutator for a unitary operator, commutators of higher order and the properties of its spectral measure, which is synthesized by Theorem \ref{thm1} and Corollary \ref{cor} below.  These are the unitary counterparts of \cite{jmp} Theorems 2.2. and 4.2. We underline that the ultimate improvments obtained in \cite{abc} are due to a systematic exploitation of the unitary framework, framework that has been logically pushed forward in the present manuscript.

The former abstract results are applied to three models. First, we propose an operator-theoretic approach to estimate the decay of the correlation functions of the Bernouilli shifts. Then, we obtain some complementary results concerning local perturbations of the Floquet operator associated to a quantum harmonic oscillator under a resonant AC-Stark potential. Third, we study the spectral properties of some GGT matrices with asymptotically constant Verblunsky coefficients, complementing various results scattered throughout the literature.

The manuscript is structured as follows. The abstract results are presented in Section 2. Sections 3, 4 and 5 are dedicated to the applications mentioned previously. The proof of Theorem \ref{thm1} is developped in Section 6.  Some auxiliary results and technicalities have been postponed in Section 7.\\

\noindent {\bf Notations:} Let us fix some notations adopted throughout this paper. Our unitary operator is defined on some fixed infinite-dimensional Hilbert space ${\cal H}$ on ${\mathbb C}$. The resolvent set of a closed operator $B$ on ${\cal H}$ is denoted by $\rho(B)$ and its spectrum by: $\sigma(B) \equiv {\mathbb C}\setminus \rho(B)$. The open unit disk and the unit circle are denoted by ${\mathbb D}$ and $\partial {\mathbb D}={\mathbb S}$ respectively. The one-dimensional torus is denoted by ${\mathbb T}$. The positive constants independent of the relevant parameters of the problem are generically denoted by $c$ or $C$. If $A$ is a self-adjoint operator defined on ${\cal H}$ with domain ${\cal D}(A)$, we use the japanese bracket notation: $\bra A\ket = \sqrt{(A^2+1)}$. Lastly, for any function $\Phi$ on ${\mathbb S}$ is associated in a unique manner to the function $\phi$ defined on ${\mathbb T}$ by:  $\phi(\theta)=\Phi(e^{i\theta})$, for all $\theta \in {\mathbb T}$. If $U$ is a unitary operator defined on ${\cal H}$ and if its spectral family is denoted by $(E_{\Delta})_{\Delta \in {\cal B}({\mathbb T})}$, where ${\cal B}({\mathbb T})$ stands for the family of Borel sets of ${\mathbb T}$, we will have that:
\begin{equation*}
\Phi(U)=\int_{\mathbb T}\phi(\theta)dE(\theta)=\int_{\mathbb T}\Phi(e^{i\theta})dE(\theta)\enspace .
\end{equation*}
We will identify frequently the spectrum of $U$ and its component (which are subsets of ${\mathbb S}$) with the corresponding support of the spectral measure, which lies in ${\mathbb T}$.

\section{Hypotheses and Main Results}

In this section, we introduce the main abstract result of this manuscript i.e. Theorem \ref{thm1}. The core of its development relies on the existence of a self-adjoint operator $A$, densely defined on ${\cal H}$ (the conjugate operator), which respect to which our unitary operator $U$ satisfies some suitable regularity conditions. We start by describing them.

\begin{defin}\label{c1} Let $B\in {\cal B}({\cal H})$ and $A$ a self-adjoint operator defined on ${\cal H}$ with domain ${\cal D}(A)$. The operator $B$ is of class $C^1$ with respect to $A$ (or shortly $B\in C^1(A)$), if there exists a dense linear subspace ${\cal S}\subset {\cal H}$, such that ${\cal S}\subset {\cal D}(A)$ and the sesquilinear form $F$ defined by
\[F(\varphi,\phi):=\langle A\varphi,B\phi\rangle - \langle \varphi,BA\phi\rangle\]
for any $(\varphi,\phi)\in {\cal S}\times {\cal S}$, extends continuously to a bounded form on ${\cal H}\times {\cal H}$. The bounded linear operator associated to the extension of $F$ is denoted by $\mathrm{ad}_A (B) =[A,B]$.
\end{defin}

\begin{defin}\label{ck} Let $k\in {\mathbb N}$, $B\in {\cal B}({\cal H})$ and $A$ a self-adjoint operator defined on ${\cal H}$ with domain ${\cal D}(A)$. The operator $B$ is of class $C^k$ with respect to $A$ (or shortly $B\in C^k(A)$), if there exists a dense linear subspace ${\cal S}\subset {\cal H}$ such that ${\cal S}\subset {\cal D}(A)$ and:
\begin{itemize}
\item $B\in C^{k-1}(A)$
\item the sesquilinear form $F$, defined by: $F(\varphi,\phi):=\langle A\varphi, \mathrm{ad}_A^{k-1}(B) \phi\rangle - \langle \varphi, \mathrm{ad}_A^{k-1}(B) A\phi\rangle $, for any $(\varphi,\phi)\in {\cal S}\times {\cal S}$, extends continuously to a bounded form on ${\cal H}\times {\cal H}$.
\end{itemize}
The bounded linear operator associated to the extension of $F$ is denoted by $\mathrm{ad}_A (\mathrm{ad}_A^{k-1}(B)) =\mathrm{ad}_A^{k}(B)$. If $B$ belongs to $C^k(A)$ for any $k\in {\mathbb N}$, we say that $B\in C^{\infty}(A)$. 
\end{defin}

Actually, the notation takes its origin in the fact that a bounded linear operator $B$ belongs to $C^k(A)$ if and only if the strongly continuous application $t \mapsto e^{itA}B e^{-itA}$ with values in ${\cal B}({\cal H})$ is strongly $C^k$ on ${\mathbb R}$. We refer to Section 7 or \cite{abmg} for more details. We shall write naturally: $C^0(A)={\cal B}({\cal H})$ and $\mathrm{ad}_A^0 B=B$. ${\cal S}$ can be equivalently chosen as ${\cal D}(A)$ in Definitions \ref{c1} and \ref{ck}. Some properties of the classes $C^k(A)$ are summed up in Section 7. 

If $U$ is a unitary operator defined on some Hilbert space ${\cal H}$, $U\in C^1(A)$ if and only if $U^*\in C^1(A)$. In particular, $U({\cal D}(A))$ and $U^*({\cal D}(A))$ are subsets of ${\cal D}(A)$, which implies that: $U({\cal D}(A))=U^*({\cal D}(A))={\cal D}(A)$. These considerations motivates the following equivalence, proved in \cite{abc} Section 5:
\begin{lem}\label{equiv1} Let $U$ be a unitary operator defined on ${\cal H}$. Then, the following assertions are equivalent:
\begin{itemize}
\item[(a)] $U\in C^1(A)$.
\item[(b)] $U^*\in C^1(A)$.
\item[(c)] There exists a dense linear subspace ${\cal S}_1$ of ${\cal H}$ such that $U {\cal S}_1= {\cal S}_1$, ${\cal S}_1 \subset {\cal D}(A)$ and the sesquilinear form $F_1:{\cal S}_1\times {\cal S}_1\rightarrow {\mathbb C}$: $F_1(\varphi,\phi):=\langle U\varphi,AU\phi\rangle\ - \langle\varphi,A\phi\rangle$ extends continuously to a bounded form on ${\cal H}\times {\cal H}$. This extension is associated to a bounded operator denoted by $U^{*}AU-A$.
\item[(d)] There exists a dense linear subspace ${\cal S}_2$ of ${\cal H}$ such that $U {\cal S}_2= {\cal S}_2$, ${\cal S}_2 \subset {\cal D}(A)$ and the sesquilinear form $F_2:{\cal S}_2\times {\cal S}_2\rightarrow {\mathbb C}$: $F_2(\varphi,\phi):=\langle \varphi,A\phi\rangle - \langle U^*\varphi,AU^*\phi\rangle$ extends continuously to a bounded form on ${\cal H}\times {\cal H}$. This extension is associated to a bounded operator denoted by $A-UAU^*$.
\end{itemize}
Moreover, $U^{*}AU-A = U^* (\mathrm{ad}_A U)$, $(\mathrm{ad}_A U) U^*=A-UAU^*$.
\end{lem}

When speaking about the positivity conditions we are about to introduce, we will write indifferently $U^{*}AU-A$ for $U^* (\mathrm{ad}_A U)$ and $(\mathrm{ad}_A U) U^*$ for $A-UAU^*$ (in the sense of Lemma \ref{equiv1}):
\begin{defin}\label{propagating} Let $A$ be a self-adjoint operator with domain ${\cal D}(A) \subset {\cal H}$ and $U$ a unitary operator which belongs to $C^1(A)$. Then, we say that given $\Theta \in {\cal B}({\mathbb T}),$
\begin{itemize}
\item {\bf $P_w$}: $U$ is weakly propagating with respect to $A$ if $U^{*}AU-A > 0$ (i.e non-negative and injective) 
\item {\bf $P({\Theta})$}: $U$ is propagating with respect to the observable $A$ on ${\Theta}$ or on the arc $e^{i\Theta}$ if there exist $c >0$ and a compact operator $K$ such that: $E_{\Theta} (U^{*}AU-A) E_{\Theta} \geq c E_{\Theta} + K$
\item {\bf $P_s(\Theta)$}: $U$ is strictly propagating with respect to the observable $A$ on ${\Theta}$ or on the arc $e^{i\Theta}$ if there exist $c >0$ such that: $E_{\Theta} (U^{*}AU-A) E_{\Theta} \geq c E_{\Theta}$.
\end{itemize}
\end{defin}
We have clearly that: ${\bf P_s({\mathbb T})}\Rightarrow {\bf P_w}$. Sometimes, we write that the operator $U$ is (strictly) propagating for $A$ at a point $\theta$ of the torus ${\mathbb T}$, when there exists an open neighbourhood $\Theta_{\theta}$ of $\theta$ such that $U$ is (strictly) propagating for $A$ on ${\Theta}_{\theta}$. Following \cite{abcf}, this is equivalent to claim that there exist a smoothed characteristic function $\phi$ supported in $\Theta_{\theta}$, which takes value 1 on a neighbourhood of $\theta$ and a positive constant $c$ such that: 
\begin{equation*}
\Phi(U)\left(U^*AU-A\right)\Phi(U) \geq c \Phi(U)^2\enspace.
\end{equation*}

\noindent{\bf Remark:} Since the spectral projectors associated to $U$ commute with $U$ and $U^*$, the positivity conditions presented in Definition \ref{propagating} can be equivalently described writing $A-UAU^*$ in place of $U^*AU-A$. This remark will be used without any further comment.

We can formulate now a first spectral result:
\begin{thm}\label{thm2} Assume that the unitary operator $U$ is weakly propagating with respect to the self-adjoint operator $A$. Then, $\sigma_{pp}(U)=\emptyset$.
\end{thm}
This result was essentially proven under a somewhat different form in \cite{abcf2}. This is a straightforward consequence of the Virial Theorem (see Paragraph 6.1). A stronger version is proposed in Section 8.

However, strenghtening the regularity hypotheses, we can derive more precise informations on the spectral properties of $U$. The following result was proven in \cite{abc} in a more general form:
\begin{thm}\label{thm} Let $\Theta$ be an open subinterval of ${\mathbb T}$. Assume that $U$ is propagating with respect to $A$ on $\Theta$ and belongs to ${\cal C}^{1,1}(A)$. Then,
\begin{itemize}
\item $U$ has a finite number of eigenvalues in $e^{i\Theta}$. Each of these eigenvalues has a finite multiplicity.
\item For any compact set $K \subset \Theta\setminus \sigma_{pp}(U),$
\begin{equation*}
\sup_{|z|\neq 1, \arg z \in K} \|\bra A \ket^{-1}(1-zU^*)^{-1}\bra A \ket^{-1} \| < \infty \enspace .
\end{equation*}
\item The spectrum of $U$ has no singular continuous component in $e^{i\Theta}$.
\end{itemize}
If $U$ is strictly propagating with respect to $A$ on $e^{i\Theta}$, then Statement (i) can be replaced by: $U$ is purely absolutely continuous on $e^{i\Theta}$.
\end{thm}
The control of the point spectrum was again obtained by a suitable version of the Virial Theorem and only uses the fact that $U$ is propagating with respect to $A$ on $\Theta$ (see Paragraph 6.1).

These results can be extended as follows, in complete analogy with the developments of the theory for self-adjoint operators. If the regularity of the unitary operator $U$ is better than expected, the conclusions of Theorem \ref{thm} are strengthened. The following result is the counterpart of \cite{jmp} Theorem 2.2:
\begin{thm}\label{thm1} Let $\Theta$ be an open subset of ${\mathbb T}$. Assume $U$ is propagating with respect to $A$ on $e^{i\Theta}$ and that there exists $k\in {\mathbb N}$ such that $U \in C^{k+1}(A)$. Let $s>k+1/2$ and $\phi \in C_0^{\infty}(\Theta \setminus \sigma_{pp}(U))$. Then, 
\begin{itemize}
\item[(i)] The conclusions of Theorem \ref{thm} hold.
\item[(ii)] For any compact subset $K\subset \Theta\setminus \sigma_{pp}(U),$
\begin{equation*}
\sup_{|z|\neq 1, \arg z \in K} \|\bra A \ket^{-s}(1-zU^*)^{-1}\bra A \ket^{-s} \| < \infty \enspace .
\end{equation*}
\item[(iii)] If $z$ tends to $e^{i\theta}$, then $\bra A \ket^{-s}(1-zU^*)^{-1}\bra A \ket^{-s}$ converges in norm to a bounded operator denoted $F^+_{1,s}(0^+,e^{i\theta})$ (resp. $F^-_{1,s}(0^+,e^{i\theta}) )$ if $|z|<1$ (resp. $|z|>1$). This convergence is uniform if $\theta$ belongs to any compact subset $K\subset \Theta\setminus \sigma_{pp}(U)$.
\item[(iv)] The operator-valued functions defined by $F^{\pm}_{1,s}$ are of class $C^k$ on each connected component of $\Theta\setminus \sigma_{pp}(U)$, with respect to the norm topology on ${\cal B}({\cal H})$.
\item[(v)] there exists $C>0$ such that for all $m\in {\mathbb Z},$
\begin{equation*}
\|\bra A \ket^{-s} U^m \Phi(U) \bra A \ket^{-s} \|\leq C \bra m \ket^{-k}\enspace.
\end{equation*}
\end{itemize}
\end{thm}
Our proof of Theorem \ref{thm1} (see Paragraph 6.4) is intrinsically based on the unitary functional calculus and is a natural extension of the developments lead in \cite{abc}. The next corollary follows from Theorem \ref{thm1} by standard interpolation arguments \cite{jmp}:
\begin{cor}\label{cor} Let $\Theta$ be an open subset of ${\mathbb T}$ and $\Phi \in C_0^{\infty}(\Theta \setminus \sigma_{pp}(U))$. Assume $U$ is propagating with respect to $A$ on $e^{i\Theta}$ and $U \in C^{\infty}(A)$. Then, for any $0<s'<s$, there exists $C>0$ such that for any $m\in {\mathbb Z}$:
\begin{equation*}
\|\bra A \ket^{-s} U^m \Phi(U) \bra A \ket^{-s} \|\leq C \bra m \ket^{-s'}\enspace.
\end{equation*}
\end{cor}

The next three sections are dedicated to the examples. The proof of Theorem \ref{thm1} is postponed to Section 6.

\section{Correlations for the Bernouilli Shifts}

In this section, we show how Corollary \ref{cor} can be reinterpreted as an operator-theoretic way to derive estimates on the correlation functions for some ergodic classical dynamical systems (see e.g. \cite{baladi} and references therein). The following development was partly borrowed from \cite{such2}.

For illustrative purposes, we have focused our discussion on a specific example, although the approach can undoubtedly be extended to any (ergodic) dynamical system for which a conjugate operator can be identified. We refer to \cite{N1}, \cite{wa} for general considerations on ergodic dynamical systems. Let $(\Omega, {\cal F}, P)$ be the following probability space: $\Omega  = \prod_{n \in {\mathbb Z}} \{-1, 1\}$, ${\cal F}$ is the $\sigma$-algebra generated by the cylinders on $\Omega$ and $P$ is the product measure, $P=\otimes_{n\in {\mathbb Z}} P_0$ where $P_0$ is the following non-trivial Bernouilli measure: for $(p,q)\in (0,1)^2,$
\begin{equation*}
P_0 = p\, \delta_{-1} + q\, \delta_{1} \quad \mathrm{with}\quad p+q=1 \enspace .
\end{equation*}
We consider also the shift $S$ on $\Omega$ defined by: $S(\omega) =  \omega'$ where $\omega'_n = \omega_{n+1}$ for all $n\in {\mathbb Z}$. Since $S$ is a measure-preserving automorphism, it is associated to a unitary Koopman operator defined on the Hilbert space $L^2(\Omega, {\cal F}, P)$ by:
\begin{equation*}
Uf = f \circ S \enspace.
\end{equation*}
The space $L^2(\Omega, {\cal F}, P)$ is an avatar of the Toy Fock space \cite{pa}, \cite{meyer}, \cite{ja}. The spectral properties of the Koopman operator $U$ are well-known: the point spectrum of $U$ is reduced to one simple eigenvalue $\{1\}$, its singular continuous component is empty and its absolutely continuous part covers the whole unit circle ${\mathbb S}$. If $Q$ denotes the orthogonal projection on the subspace generated by the constant functions on $\Omega$ and $Q^{\perp}=I-Q$, we have that: $QU=UQ=Q$ and $U=Q+Q^{\perp}UQ^{\perp}$.

We describe now a standard procedure to construct an orthonormal basis on $L^2(\Omega, {\cal F}, P)$ (the Fourier-Walsh basis). Let us choose first an orthonormal basis on $L^2(\Omega_0, 2^{\Omega_0}, P_0)\sim {\mathbb C}^2$, denoted by: $\{e_0, e_0^{\perp}\}$, where $e_0(\omega)=1$ for all $\omega \in \Omega_0$. Given any finite non-empty set $\sigma \subset {\mathbb Z}$, define on $(\Omega, {\cal F}, P)$: $f_{\sigma} = \otimes_{n \in {\mathbb Z}} f_n$ where
$$
f_n = \begin{cases} e_0, & \text{if  $ n \notin \sigma $ }  \\
                           e_0^{\perp} & \text{ if  $ n \in \sigma $ }
              \end{cases}
$$
By convention, $f_{\emptyset} =\otimes_{n \in {\mathbb Z}} e_0 \equiv 1$. The function $f_{\emptyset}$ is sometimes called the vacuum state. In particular, for $\sigma = \{i_1, \ldots, i_n\}$ and $\omega = (\omega_n)_{n=-\infty}^{\infty} \in \Omega $  one has that: $f_{\sigma}(\omega) = \prod_ {l=1}^n f_{i_l}(\omega_{i_l} )$. By construction,
$\langle f_{\sigma} , f_{\rho} \rangle = \delta_{\sigma \rho}$. If $\mathfrak{B}_n $ denotes the orthonormal set, $\mathfrak{B}_n = \{ f_{\sigma} : \sigma \subset {\mathbb Z} , |\sigma | = n \}$, ${\mathcal H}_n$ the subspace generated by $\mathfrak{B}_n $ and $Q_n$ the orthogonal projection on ${\mathcal H}_n$ ($Q_0=Q$), we have that: ${\mathcal H}_n  \perp {\mathcal H}_m $ for $m\neq n$ and
\begin{equation*}
L^2(\Omega,{\cal F}, P) = \oplus_{n=0}^{\infty} {\mathcal H}_n \enspace .
\end{equation*}
We refer again to \cite{meyer}, \cite{ja} for more details. If we denote by $L_{\perp}^2(\Omega, {\cal F}, P)$ the orthocomplement of the constant functions in $L^2(\Omega,{\cal F}, P)$, we have that: $L_{\perp}^2(\Omega, {\cal F}, P) = \{f\in L^2(\Omega, {\cal F}, P); \int_{\Omega}f\,dP= 0\} = \oplus_{n=1}^{\infty} {\mathcal H}_n$.
The operators $U$ and $U^{\perp}:=Q^{\perp}UQ^{\perp}$ are unitaries on $L^2(\Omega,{\cal F}, P)$ and $L_{\perp}^2(\Omega,{\cal F}, P)$ respectively. The action of $U$ and $U^{\perp}$ on the orthonormal basis described above is given by $Uf_{\sigma} = f_{\sigma'}$ where $\sigma' = \{ i_1+1 , \ldots, i_n+1  \}$ if $\sigma = \{i_1, \ldots, i_n\}$. The fact that $U$ and $U^{\perp}$ leave each subspace ${\mathcal H}_n$ invariant can be used in the construction of the commutation relationships as follows. For any non-negative integral number $n$, define the linear operator $A_n : {\mathcal H}_n \rightarrow {\mathcal H}_n $ as follows: $A_0 f_{\emptyset}=0$ and for any $f_{\sigma} \in {\mathcal H}_n$, $n\in {\mathbb N},$
\begin{equation*}
A_n f_{\sigma} = \left( \frac{1}{|\sigma |} \sum_{i\in \sigma} i  \right) f_{\sigma} = \left( \frac{1}{n} \sum_{i\in \sigma} i  \right) f_{\sigma}
\end{equation*}
Each operator $A_n$ is essentially self-adjoint on $\bra f_{\sigma}; |\sigma| =n\ket \subset {\mathcal H}_n$. We also denote by $A_n$ its self-adjoint extension. Due to Stone's Theorem, $(e^{itA_n})_{t\in {\mathbb R}}$ defines a strongly continuous unitary group on ${\cal H}_n$, inducing naturally two strongly continuous unitary groups on $L^2(\Omega, {\cal F}, P)$ and $L_{\perp}^2(\Omega,{\cal F}, P)$ respectively, denoted by $(\Gamma_t)_{t\in {\mathbb R}}$ and $(\Gamma_t^{\perp})_{t\in {\mathbb R}}$  and defined by:
\begin{eqnarray*}
\Gamma_t &=& \sum_{n=0}^{\infty} Q_n e^{itA_n}Q_n = \sum_{n=1}^{\infty} Q_n e^{itA_n}Q_n +Q\\
\Gamma_t^{\perp} &=& \sum_{n=1}^{\infty} Q_n e^{itA_n}Q_n \enspace .
\end{eqnarray*}
Their respective generator $A$ and $A^{\perp}$, somewhat written informally,
\begin{equation*}
A = \sum_{n=1}^{\infty} Q_n A_n Q_n = A^{\perp} \enspace ,
\end{equation*}
are self-adjoint operators on $L^2(\Omega, {\cal F}, P)$ and $L_{\perp}^2(\Omega,{\cal F}, P)$ respectively. A straightfoward computation gives: $U^*A U - A = Q^{\perp}$ on $L^2(\Omega, {\cal F}, P)$ and $U^{\perp *}A^{\perp} U^{\perp} - A^{\perp} = I$ on $L_{\perp}^2(\Omega,{\cal F}, P)$. In particular, $U^{\perp}$ is strictly propagating with respect to $A^{\perp}$ on ${\mathbb T}$ and $U^{\perp} \in C^{\infty}(A^{\perp})$. Using the fact that $Q^{\perp}$ commutes with $U$ and $A$, Corollary \ref{cor} rewrites: for any $0<s'<s$, there exists $C>0$ such that for any $m\in {\mathbb Z},$
\begin{equation*}
\|Q^{\perp}\bra A \ket^{-s} U^m \bra A \ket^{-s}Q^{\perp} \|\leq C \bra m \ket^{-s'}\enspace .
\end{equation*}

This result is a complement to the existing literature (see \cite{baladi} Chapter 1). We know that, for at least locally constant functions, this decay is exponential. Since the maps $t\mapsto e^{iAt}Ue^{-iAt}$ and $t\mapsto e^{iA^{\perp}t}U^{\perp}e^{-iA^{\perp}t}$ are norm analytic, fact that has not been exploited here, the development of an analytic version of Theorem \ref{thm1} may be appropriate to fill the gap between both approaches.

\section{Resonant AC-Stark perturbations}

Time-dependent perturbations of the quantum harmonic oscillator have been a regular subject of interest \cite{hls}, \cite{ev}, \cite{comb}, \cite{gy}, \cite{wmw}. When it is submitted to an AC-Stark potential, the Floquet operator of the system, which is explicit, undergoes a spectral transition between the resonant and non-resonant regimes \cite{ev}. The stability of these spectral properties under perturbations has been studied in the non-resonant regime \cite{wmw} and partly in the resonant regime \cite{gy}, \cite{abc}. The results obtained in this paragraph complete those obtained in \cite{abc}.

We recall the main features of the model briefly. For more details, we refer the reader to \cite{ev}. The Hamiltonian of the Harmonic oscillator (with unit mass) is defined on $L^2({\mathbb R})$ by:
$$
H_{\omega} = \frac{p^2}{2} + \frac{1}{2}\, \omega^2 x^2
$$
where $p = -i \partial_x$. We also write $\omega_0 T = 2\pi$. Let $E$ be a real-valued continuous periodic function with period $T$, $T>0$ and $(U_0 (t,s))_{(s,t)\in {\mathbb R}^2}$ be the unitary propagator associated to the AC-Stark Hamiltonian $H_0 (t) = H_{\omega} + E(t)\, x$. The propagator is explicit. In particular, for all $t\in {\mathbb R},$
\begin{equation}
U_0(t,0) = e^{-i \varphi_1(t) x}e^{i \varphi_2(t) p/\omega}e^{-iH_{\omega}t -i\psi(t)} \label{free}
\end{equation}
where $\varphi_1(t) = \int_0^t E(\tau) \cos(\omega (\tau - t)) \, d\tau$, $\varphi_2(t) =  -\int_0^t E(\tau) \sin(\omega (\tau - t ))\, d\tau$ and $\psi(t) =-\frac{1}{2}\int_0^t (\varphi_1(\tau)^2-\varphi_2(\tau)^2)\, d\tau$. For simplicity, we denote $U_0(t,0)$ by $U_0(t)$ in the following. The evolution of the observables $x$ and $p$ under this propagator are also explicit (understood on a suitable domain like the space of the Schwartz functions ${\cal S}({\mathbb R})$):
\begin{eqnarray}\label{man2}
U_0(t)^* p \, U_0(t) & = &  -x \omega \sin(\omega t) + p \cos(\omega t) + \varphi_1(t) \label{freeEnss}\\
U_0(t)^* x \, U_0(t) & = &  x \cos(\omega t) + \frac{p}{\omega} \sin(\omega t) - \frac{1}{\omega}  \varphi_2(t)\nonumber
\end{eqnarray}
These identities allows us to deduce the spectral properties of the Floquet operator $U_0(T)$: 
\begin{itemize}
\item[(1)] If $\omega_0 \neq \omega$, the Floquet operator $U_0(T)$ is pure point.
\item[(2)] If $\omega_0 =\omega$ (case coined as resonant)
\begin{itemize}
\item[(2.a)] If $\varphi_1 (T)=\varphi_2 (T)=0$, then $U_0(T)$ is pure point.
\item[(2.b)] If either $\varphi_1 (T)\neq 0$ or $\varphi_2 (T)\neq 0$, then $U_0(T)$ has purely absolutely continuous spectrum and $\sigma(U_0(T))={\mathbb S}$. Specifically,
\end{itemize}
\end{itemize}
In the case (2.a), denoting $A_1:=\varphi_1 (T)^{-1}p$ if $\varphi_1 (T)\neq 0$, $A_2:=-\omega \varphi_2 (T)^{-1}x$ if $\varphi_2 (T)\neq 0$ and ${\cal S}({\mathbb R})$ the space of the Schwarz function on ${\mathbb R}$ (observe that $U_0(T){\cal S}({\mathbb R})={\cal S}({\mathbb R})$), we have that:
\begin{itemize}
\item $U_0(T)$ belongs to $C^{\infty}(A_1)\cap C^{\infty}(A_2)$ i.e. to $C^{\infty}(p)\cap C^{\infty}(x)$.
\item If $\varphi_1 (T)\neq 0$, then the operator $U_0(T)^* A_1 U_0(T) - A_1$ defined via its sesquilinear form on ${\cal S}({\mathbb R})\times {\cal S}({\mathbb R})$ can be extended uniquely as a bounded operator on $L^2({\mathbb R})$ and $U_0(T)^* A_1 U_0(T) - A_1 = I$.
\item If $\varphi_2 (T)\neq 0$, then the operator $U_0(T)^* A_2 U_0(T) - A_2$ defined via its sesquilinear form on ${\cal S}({\mathbb R})\times {\cal S}({\mathbb R})$ can be extended uniquely as a bounded operator on $L^2({\mathbb R})$ and $U_0(T)^* A_2 U_0(T) - A_2 = I$.
\end{itemize}
We refer to \cite{abc} for the details.

In this resonant regime, some of the spectral property are preserved if the Hamiltonian $H_0(\cdot)$ is suitably perturbed. Let $V$ denote the multiplication operator by the real-valued function $V(\cdot)$ on $L^2({\mathbb R})$, $V(\cdot) \in L^{\infty}({\mathbb R})$ and define the perturbed time-dependent Hamiltonian, $H(\cdot)$ by: $H(t)= H_0(t) + V$. If the propagator $(U(t,s))$ associated to $H(\cdot)$ exists, then:
\begin{thm}\label{b} Let $\omega_0 =\omega$, $n \geq 2$ and assume that $\varphi_j (T)\neq 0$ for some $j\in \{1,2\}$. Given $V(\cdot) \in C^n({\mathbb R})$, consider the Floquet operator $U(T)$ defined by (\ref{dyson1}). Then,
\begin{itemize}
\item[(a)] If $\partial_x V(x)$ vanishes when $|x|$ tends to infinity, then there is no singular continuous component in the spectrum of $U(T)$. Moreover, its point subspace has finite dimension.
\item[(b)] if $T \|\partial_x V(\cdot)\|_{\infty} < |\varphi_1 (T)|$ (resp. $2\pi \|\partial_x V(\cdot)\|_{\infty} < |\varphi_2 (T)|$), then the spectrum of $U(T)$ is purely absolutely continuous.
\end{itemize}
If $s>n-1/2$ and and $\phi \in C_0^{\infty}(\sigma(U) \setminus \sigma_{pp}(U))$, we have in addition that:
\begin{itemize}
\item For any compact subset $K\subset \sigma (U(T)) \setminus \sigma_{pp}(U(T)),$
\begin{equation*}
\sup_{|z|\neq 1, \arg z \in K} \|\bra A_j \ket^{-s}(1-zU^*)^{-1}\bra A_j \ket^{-s} \| < \infty \enspace .
\end{equation*}
\item If $z$ tends to $e^{i\theta}$, then $\bra A_j \ket^{-s}(1-zU^*)^{-1}\bra A_j \ket^{-s}$ converges in norm to a bounded operator denoted $F^+_{1,s}(0^+,e^{i\theta})$ (resp. $F^-_{1,s}(0^+,e^{i\theta}) )$ if $|z|<1$ (resp. $|z|>1$). This convergence is uniform if $\theta$ belongs to any compact subset $K\subset \sigma (U(T)) \setminus \sigma_{pp}(U)$.
\item The operator-valued functions defined by $F^{\pm}_{1,s}$ are of class $C^{n-1}$ on each connected component of $\sigma (U(T))\setminus \sigma_{pp}(U)$, with respect to the norm topology on ${\cal B}({\cal H})$.
\item there exists $C>0$ such that for all $m\in {\mathbb Z},$
\begin{equation*}
\|\bra A_j \ket^{-s} U^m \Phi(U) \bra A_j \ket^{-s} \|\leq C \bra m \ket^{-n+1}\enspace.
\end{equation*}
\end{itemize}
\end{thm}
Statement (b) was proven in \cite{gy} for smooth and mildly unbounded potentials $V$. Statement (a) was established under weaker hypotheses in \cite{abc}.

In the proof of Theorem \ref{b}, we will restrict our discussion to the case $\varphi_1(T)\neq 0$. The other case can be treated similarly. If the propagator $(U(t,s))$ associated to $H(\cdot)$ exists, then it satisfies for all $(s,t)\in {\mathbb R}^2$, $U(t,s) = U_0(t,s) \Omega(t,s)$ where $\Omega (t,s)$ is defined in the strong sense by:
\begin{equation}\label{dyson1}
\Omega(t,s)-I = -i \int_s^t U_0(\tau,s) V U_0^*(\tau,s) \Omega(\tau,s)\, d\tau \enspace .
\end{equation}
We will denote $\Omega(t):=\Omega(t,0)$ and $U(t):=U(t,0)$. The properties of the function $V$ are related to the regularity of $U(T)$ with respect to $A_1$ (i.e. $p$) and some compactness properties as follows:
\begin{prop}\label{lemEnss0} Let $n \in {\mathbb N}$ and $V(\cdot)$ be a real-valued function in $C^n({\mathbb R})$ such that for all $k\in \{0,\ldots,n\}$, $\partial_x^k V(\cdot) \in L^\infty (\mathbb{R})$. Then,
\begin{itemize}
\item for all $t\in {\mathbb R}$, $\Omega(t)$ belongs to $C^n(p)$. In particular, for all $t\in {\mathbb R}$, $\Omega^*(t) p \Omega(t) -p$ is bounded and:
\begin{equation}
\Omega^*(t) p \Omega(t) -p = - \int_0^t \cos (\omega \tau ) U^*(\tau ) (\partial_x V) U (\tau )\, d\tau \label{eq3} \enspace.
\end{equation}
\item $U(T)$ belongs to $C^n(p)$. In particular, $U^*(T) p U(T) - p$ is bounded and
\begin{equation}\label{B(T)}
U^*(T) p U(T) - p = \varphi_1(T) - \int_0^T \cos (\omega \tau ) U^* (\tau) (\partial_x V) U (\tau)\, d\tau \enspace .
\end{equation}
\item If in addition, $\lim_{|x| \to \infty} \partial_x V(x) = 0$, the bounded operator $\Omega^*(t) p \Omega(t) -p$ is compact for any $t\in {\mathbb R}$.
\end{itemize}
\end{prop}
We refer to \cite{abc} for the proof.

\noindent{\bf Proof of Theorem \ref{b}:} It follows from the hypotheses and Proposition \ref{lemEnss0} that:
\begin{itemize}
\item $U(T)$ is propagating with respect to $A_1$ on ${\mathbb T}$ if
\begin{equation*}
\lim_{|x| \to \infty} \partial_x V(x) = 0 \enspace .
\end{equation*}
\item $U(T)$ is strictly propagating with respect to $A_1$ on ${\mathbb T}$ if $T \|\partial_x V(\cdot)\|_{\infty} < |\varphi_1 (T)|$.
\end{itemize}
We also have that $U(T)\in C^n(A_1)$. The proof follows from Theorem \ref{thm1}. \ep

\section{GGT Matrices}

GGT matrices appeared first in the theory of orthogonal polynomials on the unit circle \cite{g}. For an introduction to this subject in general and the model in particular, the reader is referred to \cite{S1}. The spectral analysis of such matrices has been undertaken in the contexts of periodic and random Verblunsky coefficients \cite{gt}, \cite{gnva1}, \cite{gnva2}, \cite{gol}, \cite{gn}. These developments are based on the theory of orthogonal polynomials and the associated transfer matrices formalism. In this section, we reconsider those GGT matrices with asymptotically constant Verblunsky coefficients by means of commutation relationships and complete the results obtained in \cite{gnva1}, \cite{gnva2}.

Our description of the model follows \cite{gt}. In this section, $(e_k)_{k\in {\mathbb Z}}$ denotes the canonical orthonormal basis of $l^2({\mathbb Z})$. The operators $T$ and $A$ denote respectively the shift and the position operator defined by:
\begin{eqnarray*}
T e_k &=& e_{k+1}\\
A e_k &=& k e_k
\end{eqnarray*}
for all $k\in {\mathbb Z}$. The reader will note that the shift operator $T$ belong to $C^{\infty}(A)$ and that for any nonnegative integral number $l$, $\mathrm{ad}_A^l(T)=T^l$ and $\mathrm{ad}_A^l(T^*)=(-1)^l T^{*l}$. Consider two sequences $(a_k)_{k\in {\mathbb Z}}$ and $(\alpha_k)_{k\in {\mathbb Z}}$ of positive and complex numbers respectively such that: $a_k^{-2}+|\alpha_k|^2= 1$ and
\begin{equation*}
\sum_{k=0}^{\infty} |\alpha_k|^2 = \infty = \sum_{k=-1}^{-\infty} |\alpha_k|^2 \enspace .
\end{equation*}
It follows from \cite{gt} Lemma 2.2 that the linear operator $H(\alpha)$ defined by:
\begin{eqnarray*}
H(\alpha) e_k = \frac{1}{a_k}e_{k-1} -\overline{\alpha_k} \sum_{i=k}^{\infty} \alpha_{i+1} \prod_{j=k+1}^{i} \frac{1}{a_j} e_i
\end{eqnarray*}
is unitary on $l^2({\mathbb Z})$. The operator $H(\alpha)$ is the GGT representation associated to the sequence of Verblunsky coefficients $(\alpha_k)$. In order to keep the amount of technicalities to a reasonable size, we will assume throughout this section that: $\inf_{k\in {\mathbb Z}} |\alpha_k|>0$. Under this assumption, the operator can be rewritten as follows:
\begin{equation}\label{Uzero}
H(\alpha) = T^*D_2(\alpha) - T^*D_1(\alpha)T(I-D_2(\alpha)T)^{-1}D_1(\alpha)^* \, ,
\end{equation}
where $D_1(\alpha)$ and $D_2(\alpha)$ are the bounded diagonal operators defined on $l^2({\mathbb Z})$ by: $D_1(\alpha) e_k = \alpha_k e_k$ and $D_2(\alpha) e_k = a_k^{-1} e_k$ for all $k\in {\mathbb Z}$.

When the sequence of Verblunsky coefficients is constant, say equal to $\alpha_{\infty}\notin \{0,1\}$ and $a>1$ is such that $|\alpha_{\infty}|^2 + a^{-2}= 1$, the associated GGT representation, denoted $H_a$, may be rewritten:
\begin{equation}\label{Uzeroc}
H_a= \frac{1}{a}\, T^* - |\alpha_{\infty} |^2 \, \sum_{j=0}^{\infty}  \left( \frac{T}{a}\right)^j= F_a(T)
\end{equation}
where the complex-valued function $F_a$ is defined on ${\mathbb C}\setminus \{0,a\}$ by:
\begin{equation}
F_a(z) = \frac{1}{az} - |\alpha |^2 \left(1-\frac{z}{a}\right)^{-1} = \frac{1-az}{z(a-z)}\, .
\end{equation}
In other words, $H_a={\cal F}^* f_a(\cdot){\cal F}$ where $f_a(\cdot)$ denotes the multiplication operator by the smooth function $f_a$ on $L^2({\mathbb T})$ and defined by:
\begin{equation}
f_a(\theta)=F_a(e^{i\theta}) = \frac{e^{-i\theta}-a}{a-e^{i\theta}} \enspace ,
\end{equation}
for $\theta \in {\mathbb T}$. Let us denote for all $\theta \in {\mathbb T}$, $g_a(\theta)= 2a\cos \theta -2$ and by $G_a$ the unique function defined on the unit circle such that $g_a(\theta)=G_a(e^{i\theta})$ for any $\theta \in {\mathbb T}$. The symmetric operator \begin{eqnarray*}
B_a &:=& G_a(T) A + A G_a(T)={\cal F}^*(-i(g_a(\theta)\partial_{\theta}+\partial_{\theta} g_a(\theta)) {\cal F} \\
&=& (aT+aT^*-2)A+A(aT+aT^*-2) \enspace ,
\end{eqnarray*}
defined on ${\cal D}=\bra e_k; k\in {\mathbb Z}\ket$ is essentially self-adjoint on this domain and its self-adjoint extension will also be denoted by $B_a$. Let us make a couple of additional observations:
\begin{itemize}
\item Since the multiplication operator by the smooth function $f_a$ on $L^2({\mathbb T})$ belongs clearly to $C^{\infty}(-i\partial_{\theta})$ and $C^{\infty}(-i(g_a(\theta)\partial_{\theta}+\partial_{\theta} g_a(\theta)))$, $H_a$ belongs to $C^{\infty}(A) \cap C^{\infty}(B_a)$.
\item The symbol $f_a$ being continuous, $\sigma (H_a)= \sigma_{ess} (H_a)= \mathrm{Ran} f_a=\Theta_a:=\{e^{i\theta}; \arg f_a(-\theta_a) \leq \theta \leq \arg f_a(\theta_a)\}$ where $\theta_a:=\cos^{-1} a^{-1}$.
\item The set of critical point of the function $f_a$ is reduced to the set $\{\pm \theta_a\}$. The multiplication operator by the smooth function $f_a$ and incidentally $H_a$ are purely absolutely continuous.
\end{itemize}

The last affirmation can also be derived from Proposition \ref{BB++}. We show in the next result how some local perturbations of $H_a$ through local fluctuations of the sequence of Verblunsky coefficients $(\alpha_k)$ may modify the spectral properties of the corresponding GGT representation. Let us introduce the family of seminorms $(p_{n_1,n_2})$ and $(q_n)$ defined (for non-negative integral numbers $n_1$, $n_2$ and $n$) on ${\mathbb C}^{\mathbb Z}$ by:
\begin{equation*}
p_{n_1,n_2}(u)= \sup_{k\in {\mathbb Z}} | k^{n_1} (\Delta^{n_2} u)_k |
\end{equation*}
where $(\Delta u)_k = u_k -u_{k-1}$ for all $k\in {\mathbb Z}$ and
\begin{equation}\label{seminorm}
q_n(u)= \sum_{m= 0}^n p_{m,m} \enspace .
\end{equation}

\begin{thm}\label{ggtperturbed2} Let $(\alpha_k) \in {\mathbb D}^{\mathbb Z}$ such that: $0<\inf_{k\in {\mathbb Z}} |\alpha_k|\leq \sup_{k\in {\mathbb Z}} |\alpha_k| < 1$. Assume that for all $k\in {\mathbb Z}$, $\alpha_k =\alpha_{\infty}(1+\delta_k)$ where $\lim_{|k|\rightarrow \infty}\delta_k =0$ and $q_n(\delta) < \infty$ for some $n\geq 2$. Then, we have the following:
\begin{itemize}
\item[(a)] $\sigma_{\text{ess}}(H(\alpha))=\sigma_{\text{ess}}(H_a)=\Theta_a$
\item[(b)] There is at most a finite number of eigenvalues in any compact subarc $K\subset \Theta_a$, $K \cap \partial \Theta_a=\emptyset$. Each of these eigenvalues has finite multiplicity.
\item[(c)] There is no singular continuous spectrum in $\Theta_a$ and $\sigma_{ac}(H(\alpha))=\Theta_a$.
\end{itemize}
If $s> n-1/2$ and $\Phi \in C_0^{\infty}(\Theta_a \setminus (\sigma_{pp}(H(\alpha)) \cup \partial \Theta_a))$, we have in addition that:
\begin{itemize}
\item For any compact subset $K\subset \Theta_a\setminus (\sigma_{pp}(H(\alpha)) \cup \partial \Theta_a),$
\begin{equation*}
\sup_{|z|\neq 1, \arg z \in K} \|\bra B_a \ket^{-s}(1-zH(\alpha)^*)^{-1}\bra B_a \ket^{-s} \| < \infty \enspace .
\end{equation*}
\item If $z$ tends to $e^{i\theta}$, then $\bra B_a \ket^{-s}(1-zH(\alpha)^*)^{-1}\bra B_a \ket^{-s}$ converges in norm to a bounded operator denoted $F^+_{1,s}(0^+,e^{i\theta})$ (resp. $F^-_{1,s}(0^+,e^{i\theta}) )$ if $|z|<1$ (resp. $|z|>1$). This convergence is uniform if $\theta$ belongs to any compact subset $K\subset \Theta_a\setminus (\sigma_{pp}(H(\alpha)) \cup \partial \Theta_a)$.
\item The operator-valued functions defined by $F^{\pm}_{1,s}$ are of class $C^k$ on each connected component of $\Theta_a \setminus (\sigma_{pp}(H(\alpha)) \cup \partial \Theta_a)$, with respect to the norm topology on ${\cal B}({\cal H})$.
\item there exists $C>0$ such that for all $m\in {\mathbb Z},$
\begin{equation*}
\|\bra B_a \ket^{-s} H(\alpha)^m \Phi(H(\alpha)) \bra B_a \ket^{-s} \|\leq C \bra m \ket^{-n+1}\enspace.
\end{equation*}
\end{itemize}
\end{thm}
Statements (a), (b), (c) are already known under weaker hypotheses \cite{gnva1}, \cite{gnva2}, \cite{abc2}. Here, the construction of the conjugate operator given differs slightly from that presented in \cite{abc2}. \\

\noindent{\bf Remark:} The eigenvalues may accumulate at the endpoints of the arc $\Theta_a$, from the discrete side or the essential side of the spectrum. If the perturbation is small, the location of these eigenvalues can also be controlled. \\

In the following, we relate the hypotheses of Theorem \ref{ggtperturbed2} with the framework of Theorem \ref{thm1}. This takes the form of three intermediary results:
\begin{prop}\label{BB++} Let $(\alpha,a)\in {\mathbb D}^*\times (1,\infty)$ such that: $|\alpha|^2 +a^{-2}= 1$. Denote by $H_a$ the associated GGT representation and $\theta_a:=\arccos (a^{-1})$. Then,
\begin{itemize}
\item $H_a$ is weakly propagating w.r.t $B_a$.
\item $H_a$ is strictly propagating w.r.t $B_a$ on any subarc of $\Theta_a$ which does not contain any endpoint of $\Theta_a$.
\end{itemize}
As a consequence, $H_a$ has purely absolutely continuous spectrum.
\end{prop}
\noindent{\bf Proof:} We know that $H_a\in C^1(A)\cap C^1(B_a)$. In particular, $\mathrm{ad}_A H_a =  (-aT-aT^* +2)(a-T)^{-2}$. Since $\mathrm{ad}_{B_a} H_a = G_a(T) (\mathrm{ad}_A H_a) +(\mathrm{ad}_A H_a)G_a(T)$ and $H_a^* B_a H_a -B_a= H_a^* (\mathrm{ad}_{B_a} H_a)$, we have that:
\begin{equation*}
H_a^* B_a H_a -B_a:=H_a^* [B_a, H_a] = 2(a-T)^{-1}(aT+aT^* -2)^2 (a-T^*)^{-1} \enspace .
\end{equation*}
Using the spectral representation of $T$, $T= \int_{\mathbb T} e^{i\theta} dE_T(\theta)$, it follows that for any $\varphi\in l^2({\mathbb Z}),$
\begin{equation*}
\bra \varphi, (H_a^* B_a H_a -B_a) \varphi \ket = \int_{\mathbb T} j_a(e^{i\theta}) d\mu_{T,\varphi}(\theta)
\end{equation*}
where $\mu_{T,\varphi}=\bra \varphi, E_T(\theta)\varphi \ket$, $j_a(\theta) =8(a\cos \theta-1)^2 |a-e^{i\theta}|^{-2}$ for all $\theta \in {\mathbb T}$. The measure $\mu_{T,\varphi}$ is purely absolutely continuous, the integrand is a non-negative continuous function which vanishes only on a finite subset of ${\mathbb T}$. Therefore, for any non-trivial vector $\varphi$, $\bra \varphi, (H_a^* B_a H_a -B_a) \varphi \ket>0$. The first statement follows. Let $\Theta$ be any open subarc of the unit circle such that $\overline{\Theta} \cap \partial\Theta_a=\emptyset$. We will prove that there exists $c>0$ such that:
\begin{equation*}
E_{H_a}(\Theta) \left( H_a^* B_a H_a -B_a \right) E_{H_a}(\Theta) \geq c E_{H_a}(\Theta)\enspace.
\end{equation*}
Indeed, for any $\varphi\in l^2({\mathbb Z}),$
\begin{equation*}
\bra E_{H_a}(\Theta) \varphi, (H_a^* B_a H_a -B_a) E_{H_a}(\Theta) \varphi \ket = \int_{\mathbb T} j_a(\theta) \chi_{f_a^{-1}(\Theta)}(\theta) d\mu_{T,\varphi}(\theta) \enspace .
\end{equation*}
The function $j_a$ is continuous and does not vanish on the compact set $f_a^{-1}(\overline{\Theta})$. Therefore, there exists $c>0$ such that for all $\theta \in f_a^{-1}(\overline{\Theta})\subset {\mathbb T}$, $j_a(\theta) \geq c$. In particular, we have that:
\begin{equation*}
\bra E_{H_a}(\Theta) \varphi, (H_a^* B_a H_a -B_a) E_{H_0}(\Theta) \varphi \ket \geq c \int_{\mathbb T} \chi_{f_a^{-1}(\Theta)}(\theta) d\mu_{T,\varphi}(\theta) = c \int_{\mathbb T} \chi_{\Theta}(\theta) d\mu_{H_0,\varphi}(\theta)
\end{equation*}
which implies our second affirmation. Since $H_a\in C^{\infty}(B_a)$, we can apply Theorem \ref{thm1} and deduce that $H_a$ is purely absolutely continuous on $\Theta$. Therefore, $H_a$ has no singular continuous spectrum in the arc $\{e^{i\theta}; \arg f_a(-\theta_a) \leq \theta \leq \arg f_a(\theta_a)\}$. On the other hand, it follows from Theorem \ref{thm2} that the point spectrum of $H_a$ is empty: the last statement follows. \ep

\begin{lem}\label{ad-compact} Let ${\cal H}$ be a Hilbert space and $A$ a self-adjoint operator defined on ${\cal H}$ with dense domain ${\cal D}(A)$. If $C$ is a compact operator on ${\cal H}$ which belongs to $C^2(A)$, then $\mathrm{ad}_A C$ is also compact.
\end{lem}
The proof is actually the remark (ii) made in the proof of \cite{abmg} Theorem 7.2.9. Due to the inclusions (5.2.10) noted in \cite{abmg}, $\mathrm{ad}_A C$ can be expressed as the norm-limit when $\varepsilon$ tends to 0, of the family of compact operators $(-i\varepsilon^{-1}(e^{iA\varepsilon}C e^{-iA\varepsilon}-C))_{\varepsilon >0}$.

\begin{lem}\label{prop-compact} Let ${\cal H}$ be a Hilbert space and $A$ a self-adjoint operator defined on ${\cal H}$ with dense domain ${\cal D}(A)$. Let $U$ and $V$ be two unitary operators defined on ${\cal H}$, which belong to $C^1(A)$, and such that $U-V$ and $\mathrm{ad}_A (U-V)$ are compact. Then, $U$ is propagating w.r.t $A$ on some open interval $\Theta \subset {\mathbb T}$ iff $V$ is propagating w.r.t $A$ on $\Theta$.
\end{lem}
\noindent{\bf Proof:} Due to the symmetry of the problem, it is enough to prove that $V$ is propagating on $\Theta$ if $U$ is propagating on $\Theta$. Let $\phi \in C_0^{\infty}(\Theta,[0,\infty))$. Since $(U^* AU-A) - (V^* AV-A) = (U-V)^* \mathrm{ad}_A U + V^* \mathrm{ad}_A (U-V)$, it follows from the hypotheses that the differences $(U^* AU-A) - (V^* AV-A)$ and $\Phi(U)(U^* AU-A)\Phi(U) - \Phi(U)(V^* AV-A)\Phi(U)$ are compact. On the other hand, $\Phi(U)(V^* AV-A)\Phi(U) = \Phi(V)(V^* AV-A)\Phi(V)+ (\Phi(U)-\Phi(V))(V^* AV-A)\Phi(V)+\Phi(U)(V^* AV-A)(\Phi(U)-\Phi(V))$. Since $U-V$ is compact, $\Phi(U)-\Phi(V) $ is compact (see e.g. \cite{abcf} Lemma 4.1 or Stone Weierstrass Theorem). Therefore, $\Phi(U)(U^* AU-A)\Phi(U) - \Phi(V)(V^* AV-A)\Phi(V)$ is compact. If $U$ is propagating w.r.t $A$ on $\Theta$, then there exists $c>0$ such that: $\Phi(U)(U^*AU-A)\Phi(U)\leq c \Phi(U)^2+K$ for some compact operator $K$. Note that $\Phi^2(U)-\Phi^2(V)=\Phi(U)^2-\Phi(V)^2$ is also compact since $\phi^2 \in C_0^{\infty}(\Theta,[0,\infty))$ and $U-V$ is compact. The conclusion follows combining the former observations. \ep
\\

\noindent{\bf Proof of Theorem \ref{ggtperturbed2}:} It follows from Corollary \ref{criteria1} and the hypotheses of Theorem \ref{ggtperturbed2}, that $D_1(\alpha)$ belongs to $C^n(B_a)$ and the difference $D_1(\alpha)-\alpha_{\infty}$ is compact. Due to Lemma \ref{ggt1-0}, $H(\alpha)-H_a$ is compact, meaning that the operators $H(\alpha)$ and $H_a$ have the same essential spectrum (Weyl's Theorem) and implies the first claim.
It follows from Lemmata \ref{ggt1-0}, \ref{ggt2}, \ref{ad-compact}, \ref{prop-compact} and Proposition \ref{BB++} that the unitary operator $H(\alpha)$ also belongs to $C^n(B_a)$ and is propagating for the observable $B_a$ on any open subarc $\Theta \subset \Theta_a$ such that $\overline{\Theta}\cap \partial \Theta_a=\emptyset$. The conclusion follows from Theorem \ref{thm1}. \ep

The remainder of the article is devoted to the proofs we have left aside in the previous sections. 

\section{Towards the proof of Theorem \ref{thm1}}

Although the proof of Theorem \ref{thm1} follow the lines of \cite{jmp}, it is intrinsically based on the unitary functional calculus. Its development is articulated on two axes:
\begin{itemize}
\item The control of the (embedded) point spectrum by means of the Virial Theorem (Paragraph 6.1)
\item The study of the continuous component of the spectrum using Mourre differential inequality strategy (Paragraph 6.2)
\end{itemize}
The proof is carried out in Paragraphs 6.3 and 6.4.

Before starting, let us remind or fix some notations. ${\mathbb D}^*$ will stand for ${\mathbb D}-\{ 0\}$. If $\Theta$ is an open interval in ${\mathbb T}$ and $r>1$, we denote by $S_{\Theta,r}^{\pm}$ and $\Omega_{\Theta,r}^{\pm}$ the sectors
\begin{eqnarray*}
S_{\Theta,r}^+ &=& \{ z \in {\mathbb C}; \, \arg(z) \in \Theta, r^{-1}<|z|<1\}\\
S_{\Theta,r}^- &=& \{ z \in {\mathbb C}; \, \arg(z) \in \Theta, 1<|z|<r\}\\
\Omega_{\Theta,r}^+ &=& \{ z \in {\mathbb C}; \, \arg(z) \in \Theta, r^{-1}<|z|\leq 1\}\\
\Omega_{\Theta,r}^- &=& \{ z \in {\mathbb C}; \, \arg(z) \in \Theta, 1\leq |z|<r\}\enspace.
\end{eqnarray*}
The spectral measure of $U$ is denoted by $(E(\Delta))_{\Delta \in {\cal B}({\mathbb T})}$.

Following Lemma \ref{equiv1}, the following equivalence, justified in \cite{abc} Section 5, will be used throughout this Section without any further comments:
\begin{lem}\label{equiv2} Let $k\in {\mathbb N}$ and $U$ be a unitary operator on ${\cal H}$. Then, the four following assertions are equivalent:
\begin{itemize}
\item[(a)] $U\in C^k(A)$.
\item[(b)] $U^*\in C^k(A)$.
\item[(c)] $U$ satisfies item (c) of Lemma \ref{equiv1} and $(U^* AU-A) \in C^{k-1}(A)$.
\item[(d)] $U$ satisfies item (d) of Lemma \ref{equiv1} and $(A-UAU^*) \in C^{k-1}(A)$.
\end{itemize}
\end{lem}

\subsection{The Virial Theorem and its consequences}

As mentionned at the beginning of this section, the control of the point spectrum is achieved after establishing the Virial Theorem \cite{abc}:
\begin{thm}\label{virial} Assume that $U\in C^1(A)$. Then, for all $\theta \in {\mathbb T},$, $E_{\{\theta\}} (U^* AU -A) E_{\{\theta\}} =0$. In particular, if $\varphi$ is an eigenvector of $U$, $\bra \varphi, (U^* AU -A)\varphi \ket =0$.
\end{thm}

This allows us to restate \cite{abcf} Corollary 5.1:
\begin{cor}\label{virial-2} Assume $U$ is propagating with respect to $A$ on the Borel subset $\Theta\subset {\mathbb T}$. Then, $U$ has a finite number of eigenvalues in $\Theta$. Each of these eigenvalues has finite multiplicity. 
\end{cor}
 The conclusions of Corollary \ref{virial-2} can be strenghtened under stronger hypothesis as shown in Theorem \ref{weak}.

If $U$ is propagating with respect to $A$ on some Borel subset $\Theta\subset {\mathbb T}$, it follows that $\Theta \cap \sigma_{pp}(U)$ is finite. Therefore, for any $\theta \in \Theta\setminus \sigma_{pp}(U)$, there exist $\delta_{\theta} >0$ and $c_{\theta}>0$ such that:
\begin{equation*}
E_{(\theta-2\delta_{\theta}, \theta+2\delta_{\theta})} (U^* AU -A) E_{(\theta-2\delta_{\theta}, \theta+2\delta_{\theta})} \geq c_{\theta} E_{(\theta-2\delta_{\theta}, \theta+2\delta_{\theta})} \enspace .
\end{equation*}
In other words, $U$ is strictly propagating at $\theta$. This motivates the development of the next section.

\subsection{Differential inequalities}

What follows is an adaptation of \cite{abmg} paragraph 7.3 and \cite{jmp} to our unitary formalism. From now and until the end of this paragraph, we assume that $U$ is strictly propagating with respect to some self-adjoint operator $A$ at $\theta_0 \in {\mathbb T}$:
\begin{equation*}
E_{(\theta_0 - 2\delta, \theta_0 + 2\delta)} (A-UAU^*) E_{(\theta_0 - 2\delta, \theta_0 + 2\delta)} \geq a_1 E_{(\theta_0 - 2\delta, \theta_0 + 2\delta)} \enspace ,
\end{equation*}
for some $\delta>0$, $a_1>0$ (see also Corollary \ref{equiv1}). We also assume that $(B(\varepsilon))_{\varepsilon\in (0,\varepsilon_0]}$ is a family of uniformly bounded operators on ${\cal H}$ such that: $\lim_{\varepsilon \rightarrow 0}\|B(\varepsilon)-(A-UAU^*)\| =0$.

Denoting $B_1:=A-UAU^*$, we have that:
\begin{lem}\label{lem1} There exists $C>0$ such that for all $\varepsilon \in (0,\varepsilon_0],$
\begin{gather*}
\|e^{-\varepsilon B(\varepsilon)}(e^{-\varepsilon B(\varepsilon)})^* - e^{-2\varepsilon B_1}\|\leq C\, \varepsilon \\
\|(e^{\varepsilon B(\varepsilon)})^* e^{\varepsilon B(\varepsilon)} - e^{2\varepsilon B_1}\|\leq C\, \varepsilon \enspace.
\end{gather*}
\end{lem}
Note that: $(e^{\varepsilon B(\varepsilon)})^*= e^{\varepsilon B(\varepsilon)^*}$. For $\varepsilon \in (0,\varepsilon_0]$ and $z\in \overline{\mathbb D}\setminus \{0\}$, we define:
\begin{eqnarray*}
T_{\varepsilon}^{+}(z) &=& 1- z U^* e^{-\varepsilon B(\varepsilon)}\\
T_{\varepsilon}^{-}(z) &=& 1- \bar{z}^{-1} U^* (e^{\varepsilon B(\varepsilon)})^*
\end{eqnarray*}

The following estimates are proven in \cite{abc} Section 4:
\begin{lem}\label{lem4} The linear operators $T_{\varepsilon}^{\pm}(z)$ are invertible in ${\cal B}({\cal H})$, provided $(\varepsilon,z) \in (0,\varepsilon_2]\times \Omega_{(\theta_0 - \delta, \theta_0 + \delta),2}^+$ or $(\varepsilon,z) \in [0,\varepsilon_2]\times S_{(\theta_0 - \delta, \theta_0 + \delta),2}^+$ for some $\varepsilon_2 \in (0,\varepsilon_1]$. Denote by $G_{\varepsilon}^{\pm}(z)$ the respective inverse of $T_{\varepsilon}^{\pm}(z)$. Then, there exists $C>0$, such that:
\begin{itemize}
\item For all $(\varepsilon,z) \in (0,\varepsilon_2]\times \Omega_{(\theta_0 - \delta, \theta_0 + \delta),2}^+$: $\|G_{\varepsilon}^{\pm}(z)\| \leq C \varepsilon^{-1}$.
\item For all $(\varepsilon,z) \in [0,\varepsilon_2]\times S_{(\theta_0 - \delta, \theta_0 + \delta),2}^+$: $\|G_{\varepsilon}^{\pm}(z)\| \leq C (1-|z|^2)^{-1}$.
\end{itemize}
Morever, there exists $C>0$, such that for all $(\varepsilon,z) \in (0,\varepsilon_2]\times \Omega_{(\theta_0 - \delta, \theta_0 + \delta),2}^+$ and all $\psi \in {\cal H},$
\begin{equation*}
\|G_{\varepsilon}^{\pm}(z) \psi\| \leq  C \left( \sqrt{\frac{|\bra \psi, \Re(G_{\varepsilon}^{\pm}(z)) \psi \ket |}{\varepsilon}} +\|\psi \| \right)\enspace.
\end{equation*}
\end{lem}

Let us recall two technical results:
\begin{lem}\label{lem6} Let $J\subset {\mathbb R}$ be an open bounded interval and $C$ defined by:
\begin{eqnarray*}
C: J & \rightarrow & {\cal B}({\cal H})\\
\varepsilon & \mapsto & C(\varepsilon)
\end{eqnarray*}
be a $C^1$ function with respect to the norm topology on ${\cal B}({\cal H})$. Then, the map $\varepsilon \mapsto e^{-C(\varepsilon)}$ is also norm-$C^1$ on the interval $J$. Moreover, for all $\varepsilon \in J,$
\begin{equation*}
e^{C(\varepsilon)} \partial_{\varepsilon}e^{-C(\varepsilon)}= -\sum_{p=1}^{\infty}\frac{1}{p!}\mathrm{ad}_{C(\varepsilon)}^{p-1} (\partial_{\varepsilon}{C(\varepsilon)}) \enspace .
\end{equation*}
\end{lem}
\noindent{\bf Remark:} Note that for any $\varepsilon \in (0,\varepsilon_0],$
\begin{eqnarray*}
\partial_{\varepsilon} e^{-\varepsilon B(\varepsilon)} &=& \sum_{k=1}^{\infty} \frac{(-\varepsilon)^k}{k!} \partial_{\varepsilon} (B(\varepsilon))^k - B(\varepsilon) e^{-\varepsilon B(\varepsilon)} \\
\partial_{\varepsilon} e^{\varepsilon B(\varepsilon)^*} &=& \sum_{k=1}^{\infty} \frac{\varepsilon^k}{k!} \partial_{\varepsilon} (B(\varepsilon)^*)^k - B(\varepsilon)^* e^{\varepsilon B(\varepsilon)^*} \enspace.
\end{eqnarray*}

\begin{lem}\label{8-2}{\bf Baker-Campbell-Hausdorff.} Let $k\in {\mathbb N}$. If $C \in C^k(A)$, then $e^{C} \in C^k(A)$. Moreover,
\begin{equation*}
e^{-C} A e^{C}-A=\sum_{k=1}^{\infty}  \frac{(-1)^{k-1}}{k!}\, \mathrm{ad}_C^{k-1}(\mathrm{ad}_A C) \enspace.
\end{equation*}
In particular, the following estimates hold: $\| \mathrm{ad}_A e^{C} \| \leq e^{\|C\|} \| \mathrm{ad}_A C\|$ and $\| e^{-C} A e^{C}-A \| \leq e^{\|C\|}\| \mathrm{ad}_A C \|$.
\end{lem}

The proof of Lemma \ref{8-2} is explicited in \cite{abc}. The following proposition is also proven in \cite{abc} Section 4:
\begin{prop}\label{lem7} Suppose that the map defined on $(0,\varepsilon_0]$, $\varepsilon \mapsto B(\varepsilon)$ is $C^1$ w.r.t the norm topology on ${\cal B}({\cal H})$. Then for any fixed $z\in \Omega_{(\theta_0 - \delta, \theta_0 + \delta),2}^+$, the map $\varepsilon \mapsto G_{\varepsilon}^{\pm}(z)$ is $C^1$ on $(0,\varepsilon_2]$ with respect to the norm topology. Moreover, if for any $\varepsilon \in (0,\varepsilon_0]$, $U$ and $B(\varepsilon)$ belong to $C^1(A)$ then, given $(\varepsilon,z) \in (0,\varepsilon_2]\times \Omega_{(\theta_0 - \delta, \theta_0 + \delta),2}^+$, $G_{\varepsilon}^{\pm}(z)$ belongs to $C^1(A)$ and we have that: $\partial_{\varepsilon} G_{\varepsilon}^{\pm}(z) = \pm \mathrm{ad}_A G_{\varepsilon}^{\pm}(z) + G_{\varepsilon}^{\pm}(z) Q^{\pm}(\varepsilon,z) G_{\varepsilon}^{\pm}(z)$ where,
\begin{eqnarray*}
Q^+(\varepsilon,z) &=& z U^* \left( \partial_{\varepsilon} e^{-\varepsilon B(\varepsilon)} + B_1 e^{-\varepsilon B(\varepsilon)} - \mathrm{ad}_A e^{-\varepsilon B(\varepsilon)} \right) \\
Q^-(\varepsilon,z) &=& \bar{z}^{-1} U^* \left( \partial_{\varepsilon} e^{\varepsilon B(\varepsilon)^*} -B_1 e^{\varepsilon B(\varepsilon)^*} + \mathrm{ad}_A e^{\varepsilon B(\varepsilon)^*} \right) \enspace .
\end{eqnarray*}
\end{prop}

An immediate corollary is:
\begin{cor}\label{lem7-2} Suppose that the map defined on $(0,\varepsilon_0]$, $\varepsilon \mapsto B(\varepsilon)$ is $C^1$ w.r.t the norm topology on ${\cal B}({\cal H})$ and that for any $\varepsilon \in (0,\varepsilon_0]$, $U$ and $B(\varepsilon)$ belong to $C^1(A)$. Given $(\varepsilon,z) \in (0,\varepsilon_2]\times \Omega_{(\theta_0 - \delta, \theta_0 + \delta),2}^+$, we have that for any $k\in {\mathbb N},$
\begin{eqnarray*}
\partial_{\varepsilon}(G_{\varepsilon}^{\pm}(z))^k &=& \pm \mathrm{ad}_A (G_{\varepsilon}^{\pm}(z))^k + Q_k^{\pm}(\varepsilon,z)\\
\text{where}\quad Q_k^{\pm}(\varepsilon,z) &=& \sum_{j=0}^{k-1} G_{\varepsilon}^{\pm}(z)^{j+1} Q^{\pm}(\varepsilon,z) G_{\varepsilon}^{\pm}(z)^{k-j} \enspace.
\end{eqnarray*}
\end{cor}
\noindent {\bf Proof:} Let $k\in {\mathbb N}$. In view of Lemma \ref{lem7}, we have that for any $(\varepsilon,z) \in (0,\varepsilon_2]\times \Omega_{(\theta_0 - \delta, \theta_0 + \delta),2}^+$:
\begin{equation*}
\partial_{\varepsilon}(G_{\varepsilon}^{\pm}(z))^k = \sum_{j=0}^{k-1}(G_{\varepsilon}^{\pm}(z))^j \left( \partial_{\varepsilon} G_{\varepsilon}^{\pm}(z)\right) (G_{\varepsilon}^{\pm}(z))^{k-1-j} \enspace ,
\end{equation*}
which implies the result. \ep

Now, we have all the ingredients to introduce various differential inequalities. Let us introduce more notations: given $s\in [1,\infty)$ and $k\in {\mathbb N}$, define the bounded operator-valued functions $F^{\pm}_{s,k}$ on $(0,\varepsilon_2]\times \Omega_{(\theta_0 - \delta, \theta_0 + \delta),2}^+$ by:
\begin{equation*}
F_{s,k}^{\pm}(\varepsilon,z) = \bra A\ket^{-s} (G_{\varepsilon}^{\pm}(z))^k \bra A\ket^{-s} \enspace.
\end{equation*}
\begin{lem}\label{diffin} Let $s\in [1,\infty)$. Suppose that the map defined on $(0,\varepsilon_0]$ by $\varepsilon \mapsto B(\varepsilon)$ is $C^1$ w.r.t the norm topology on ${\cal B}({\cal H})$ and that for any $\varepsilon \in (0,\varepsilon_0]$, $U$ and $B(\varepsilon)$ belong to $C^1(A)$. Then, there exists $C>0$ such that for all $(\varepsilon,z) \in (0,\varepsilon_2]\times \Omega_{(\theta_0 - \delta, \theta_0 + \delta),2}^+,$
\begin{align}\label{diffineq}
\|\partial_{\varepsilon} F_{s,1}^{\pm}(\varepsilon,z) \| \leq  C & \varepsilon q(\varepsilon) \left( \varepsilon^{-1/2} \|F_{s,1}^{\pm}(\varepsilon,z)\|^{1/2}+1 \right) \left( \varepsilon^{-1/2} \|F_{s,1}^{\mp}(\varepsilon,z)\|^{1/2}+1 \right) \nonumber \\
& + C \left( \varepsilon^{-1/2} \|F_{s,1}^{\pm}(\varepsilon,z)\|^{1/2}+ \varepsilon^{-1/2} \|F_{s,1}^{\mp}(\varepsilon,z)\|^{1/2}+1 \right)
\end{align}
with $q(\varepsilon)= \varepsilon^{-1} \max (\sup_{z \in \Omega_{(\theta_0 - \delta, \theta_0 + \delta),2}^+}\|Q^{\pm}(\varepsilon,z)\|)$.
\end{lem}
\noindent {\bf Proof:} Let $(\varphi_1,\varphi_2) \in {\cal D}(A)^2$. It follows from Lemma \ref{lem7} that:
\begin{eqnarray*}
\bra \varphi_1, \partial_{\varepsilon} G_{\varepsilon}^{\pm}(z) \varphi_2 \ket &=& \pm \bra A\varphi_1, G_{\varepsilon}^{\pm}(z) \varphi_2 \ket \mp \bra G_{\varepsilon}^{\pm}(z)^* \varphi_1, A\varphi_2 \ket  + \bra G_{\varepsilon}^{\pm}(z)^* \varphi_1, Q^{\pm}(\varepsilon,z) G_{\varepsilon}^{\pm}(z) \varphi_2 \ket \\
\mathrm{where} \quad G_{\varepsilon}^+(z)^* &=& - \bar{z}^{-1}U^* e^{\varepsilon B(\varepsilon)^*} G_{\varepsilon}^-(z) \\
G_{\varepsilon}^-(z)^* &=& - z U^* e^{\varepsilon B(\varepsilon)} G_{\varepsilon}^+(z) \enspace,
\end{eqnarray*}
for all $(\varepsilon,z) \in (0,\varepsilon_2]\times \Omega_{(\theta_0 - \delta, \theta_0 + \delta),2}^+$. We deduce that:
\begin{equation*}
|\partial_{\varepsilon} \bra \varphi_1, G_{\varepsilon}^{\pm}(z) \varphi_2 \ket | \leq C \left(\|A\varphi_1 \| \| G_{\varepsilon}^{\pm}(z) \varphi_2 \| + \|A\varphi_2 \| \| G_{\varepsilon}^{\mp}(z) \varphi_1 \| + \varepsilon q(\varepsilon) \| G_{\varepsilon}^{\pm}(z) \varphi_2 \| \| G_{\varepsilon}^{\mp}(z) \varphi_1 \| \right) \enspace.
\end{equation*}
Therefore, there exists $C>0$ such that for all $(\varepsilon,z) \in (0,\varepsilon_2]\times \Omega_{(\theta_0 - \delta, \theta_0 + \delta),2}^+$ and all $(\psi_1,\psi_2) \in {\cal H}\times {\cal H},$
\begin{equation*}
|\partial_{\varepsilon} \bra \psi_1, F_{s,1}^{\pm}(\varepsilon,z) \psi_2 \ket | \leq C \left( \| G_{\varepsilon}^{\pm}(z) \bra A\ket^{-s}\psi_2 \| + \| G_{\varepsilon}^{\mp}(z) \bra A\ket^{-s}\psi_1 \| + \varepsilon q(\varepsilon) \| G_{\varepsilon}^{\pm}(z) \bra A\ket^{-s}\psi_2 \| \| G_{\varepsilon}^{\mp}(z) \bra A\ket^{-s}\psi_1 \| \right) \enspace.
\end{equation*}
On the other hand, by Lemma \ref{lem4}, there exists $C>0$, such that for all $(\varepsilon,z) \in (0,\varepsilon_2]\times \Omega_{(\theta_0 - \delta, \theta_0 + \delta),2}^+$ and any $\psi \in {\cal H},$
\begin{equation*}
\|G_{\varepsilon}^{\pm}(z) \psi\| \leq  C \left( \sqrt{\frac{|\bra \psi, \Re (G_{\varepsilon}^{\pm}(z)) \psi \ket |}{\varepsilon}} +\|\psi \| \right)\enspace.
\end{equation*}
The conclusion follows, once noted that for any bounded operator $B$, $\|B\| =\sup_{\|\psi_1\|=1, \|\psi_2\|=1} |\bra \psi_1,B\psi_2\ket|$. \ep

We also observe that:
\begin{lem}\label{lem4-1} Let $s\in [1,\infty)$. Assume that $\sup_{(\varepsilon,z) \in (0,\varepsilon_2]\times \Omega_{(\theta_0 - \delta, \theta_0 + \delta),2}^+} \|F_{s,1}^{\pm}(\varepsilon,z) \| < \infty$. Then, there exists $C>0$ such that for all $(\varepsilon,z) \in (0,\varepsilon_2]\times \Omega_{(\theta_0 - \delta, \theta_0 + \delta),2}^+$, $\|G_{\varepsilon}^{\pm}(z) \bra A \ket^{-s}\| \leq C \varepsilon^{-1/2}$ and $\|\bra A \ket^{-s} G_{\varepsilon}^{\pm}(z)\| \leq C \varepsilon^{-1/2}$.
\end{lem}
\noindent{\bf Proof:} The first part is a straightforward consequence of Lemma \ref{lem4}. The second inequality result from the first once noted that $\|B\|=\|B^*\|$ for any $B\in {\cal B}({\cal H})$ and that:
\begin{eqnarray*}
G_{\varepsilon}^+(z)^* &=& - \bar{z}^{-1}U^* (e^{\varepsilon B(\varepsilon)})^* G_{\varepsilon}^- (z)\\
G_{\varepsilon}^- (z)^* &=& - z U^* e^{\varepsilon B(\varepsilon)} G_{\varepsilon}^+(z) \enspace .
\end{eqnarray*} \ep

The following result is similar to \cite{jmp} Theorem 2.2:
\begin{lem}\label{diffin2} Let $s\in [1,\infty)$ and $k\geq 2$. Suppose that the map defined on $(0,\varepsilon_0]$ by $\varepsilon \mapsto B(\varepsilon)$ is $C^1$ w.r.t the norm topology on ${\cal B}({\cal H})$ and that for any $\varepsilon \in (0,\varepsilon_0]$, $U$ and $B(\varepsilon)$ belong to $C^1(A)$. If
\begin{equation*}
\sup_{(\varepsilon,z) \in (0,\varepsilon_2]\times \Omega_{(\theta_0 - \delta, \theta_0 + \delta),2}^+} \|F_{1,s}^{\pm}(\varepsilon,z) \| < \infty \enspace,
\end{equation*}
then, there exists $C>0$ such that for all $(\varepsilon,z) \in (0,\varepsilon_2]\times \Omega_{(\theta_0 - \delta, \theta_0 + \delta),2}^+,$
\begin{equation}\label{diffineq-2}
\| \partial_{\varepsilon} F_{s,k}^{\pm}(\varepsilon,z) \| \leq C \left(\| F_{s,k}^{\pm}(\varepsilon,z) \|^{1-1/s} \varepsilon^{(1-2k)/2s}+q(\varepsilon) \varepsilon^{-k+1} \right)\enspace.
\end{equation}
\end{lem}
\noindent{\bf Proof:} It follows from Corollary \ref{lem7-2} that given $(\varepsilon,z) \in (0,\varepsilon_2]\times \Omega_{(\theta_0 - \delta, \theta_0 + \delta),2}^+,$
\begin{eqnarray*}
\partial_{\varepsilon} F_{s,k}^{\pm}(\varepsilon,z) &=& \bra A\ket^{-s} \mathrm{ad}_A (G_{\varepsilon}^{\pm}(z)^k) \bra A\ket^{-s} + \bra A\ket^{-s} Q_k^{\pm}(\varepsilon,z) \bra A\ket^{-s}\\
Q_k^{\pm}(\varepsilon,z) &=& \sum_{j=0}^{k-1} G_{\varepsilon}^{\pm}(z)^{j+1} Q^{\pm}(\varepsilon,z) G_{\varepsilon}^{\pm}(z)^{k-j} \enspace.
\end{eqnarray*}
In view of Lemma \ref{lem4} and Corollary \ref{lem4-1}, the second term on the RHS can be estimated by:
\begin{eqnarray*}
\|\bra A\ket^{-s} Q_k^{\pm}(\varepsilon,z) \bra A\ket^{-s}\| &\leq & k \varepsilon q(\varepsilon) \|\bra A\ket^{-s} G_{\varepsilon}^{\pm}(z)\| \| G_{\varepsilon}^{\pm}(z)\bra A\ket^{-s}\| \|G_{\varepsilon}^{\pm}(z)\|^{k-1}\\
&\leq & C_k q(\varepsilon) \varepsilon^{-k+1} \enspace ,
\end{eqnarray*}
for some $C_k>0$. Using interpolation (see e.g. \cite{fk}), Lemma \ref{lem4} and Corollary \ref{lem4-1}, we have for the first term on the RHS:
\begin{eqnarray*}
\| \bra A\ket^{-s} \mathrm{ad}_A (G_{\varepsilon}^{\pm}(z)^k) \bra A\ket^{-s}\| &\leq & \| \bra A\ket^{1-s} G_{\varepsilon}^{\pm}(z)^k \bra A\ket^{-s}\|+\| \bra A\ket^{-s} G_{\varepsilon}^{\pm}(z)^k \bra A\ket^{1-s}\| \\
&\leq & \| \bra A\ket^{-s} G_{\varepsilon}^{\pm}(z)^k \bra A\ket^{-s}\|^{1-1/s} \left(\| G_{\varepsilon}^{\pm}(z)^k \bra A\ket^{-s}\|^{1/s}+\|\bra A\ket^{-s} G_{\varepsilon}^{\pm}(z)^k\|^{1/s}\right)\\
&\leq & C \|F_{s,k}(\varepsilon,z)\|^{1-1/s} \varepsilon^{(-2k+1)/2s}\enspace,
\end{eqnarray*}
for some $C>0$, which proves the Lemma. \ep

The next step consists in integrating the differential inequalities of Lemmata \ref{diffin} and \ref{diffin2}. This is done by using an avatar of the Gronwall Lemma:
\begin{lem}\label{gronwall} Let $J=(a,b) \subset {\mathbb R}$ be an open interval and let $f$, $\varphi$ and $\psi$ be non-negative real functions on $J$ with $f$ bounded, $\varphi$ and $\psi$ in $L^1(J)$. Assume there exists $\omega \geq 0$ and $\theta \in [0,1)$ such that for all $\lambda\in J$:
\begin{equation*}
f(\lambda) \leq \omega + \int_{\lambda}^b (\varphi(\tau) f(\tau)^{\theta}+ \psi(\tau) f(\tau))\, d\tau
\end{equation*}
Then for any $\lambda\in J$, one has
\begin{equation*}
f(\lambda) \leq \left[\omega^{1-\theta}+(1-\theta) \int_{\lambda}^b \varphi(\mu) e^{(\theta-1)\int_{\mu}^b \psi(\tau)\, d\tau}\,d\mu \right]^{1/(1-\theta)} \cdot e^{\int_{\lambda}^b \psi(\tau)\, d\tau}
\end{equation*}
\end{lem}
 We refer to \cite{hm} chapter III for a proof.

As a consequence, we obtain:
\begin{lem}\label{gronwall-1} Suppose that the map defined on $(0,\varepsilon_0]$ by $\varepsilon \mapsto B(\varepsilon)$ is $C^1$ w.r.t the norm topology on ${\cal B}({\cal H})$ and that for any $\varepsilon \in (0,\varepsilon_0]$, $U$ and $B(\varepsilon)$ belong to $C^1(A)$. If
\begin{equation*}
\int_0^{\varepsilon_0} q(\varepsilon)\, d\varepsilon < \infty \enspace ,
\end{equation*}
then there exist $C>0$ and $H\in L^1((0,\varepsilon_2])$ such that for all $(\varepsilon,z) \in (0,\varepsilon_2]\times \Omega_{(\theta_0 - \delta, \theta_0 + \delta),2}^+,$
\begin{eqnarray*}
\|F_{1,1}(\varepsilon,z)\| &<& C \\
\|\partial_{\varepsilon}F_{1,1}(\varepsilon,z)\| &\leq & H(\varepsilon) \enspace.
\end{eqnarray*}
\end{lem}
\noindent {\bf Proof:} The reader will observe first that the integrability of the function $q$ implies the integrability of the function $\varepsilon \mapsto \varepsilon q(\varepsilon)$. Define, the auxiliary functions $K$ and $L$ by
\begin{eqnarray*}
K(\varepsilon,z) &=& \| F_{1,1}^+(\varepsilon,z)\| +  \| F_{1,1}^-(\varepsilon,z)\| \\
L(\varepsilon) &=& \sup_{z\in \Omega_{(\theta_0 - \delta, \theta_0 + \delta),2}^+} K(\varepsilon,z)
\end{eqnarray*}
Up some adjustment of the constants, we have that for all $(\varepsilon,z) \in (0,\varepsilon_2]\times \Omega_{(\theta_0 - \delta, \theta_0 + \delta),2}^+,$
\begin{eqnarray*}
\left| K(\varepsilon_2,z) - K(\varepsilon,z) \right| &= & \left| \| F_{1,1}^+(\varepsilon_2,z)\| - \| F_{1,1}^+(\varepsilon,z)\|+  \| F_{1,1}^-(\varepsilon_2,z)\| - \| F_{1,1}^-(\varepsilon,z)\| \right| \\
&\leq & \| F_{1,1}^+(\varepsilon_2,z) - F_{1,1}^+(\varepsilon,z)\|+  \| F_{1,1}^-(\varepsilon_2,z) - F_{1,1}^-(\varepsilon,z)\| \\
& \leq & \int_{\varepsilon}^{\varepsilon_2} \| \partial_{\rho} F_{1,1}^+(\rho,z) \|+\| \partial_{\rho} F_{1,1}^-(\rho,z) \| \, d\rho \\
& \leq & C \int_{\varepsilon}^{\varepsilon_2} (q(\rho)  K(\rho,z) + \rho^{-1/2} K(\rho,z)^{1/2} + \rho q(\rho) +1) \, d\rho 
\end{eqnarray*}
using Lemma \ref{diffin} and the fact that: $\| F_{1,1}^{\pm}(\varepsilon,z)\| \leq K(\varepsilon,z)$. It follows from Lemma \ref{lem4} that for all $(\varepsilon,z) \in (0,\varepsilon_2]\times \Omega_{(\theta_0 - \delta, \theta_0 + \delta),2}^+,$
\begin{eqnarray*}
K(\varepsilon,z) &\leq & K(\varepsilon_2,z) + C \int_{\varepsilon}^{\varepsilon_2} (q(\rho)  K(\rho,z) + \rho^{-1/2} K(\rho,z)^{1/2} + \rho q(\rho) +1) \, d\rho \\
&\leq & C\left(\varepsilon_2^{-1} +\int_{\varepsilon}^{\varepsilon_2} (q(\rho)  K(\rho,z) + \rho^{-1/2} K(\rho,z)^{1/2} + \rho q(\rho) +1) \, d\rho \right) \\
L(\varepsilon) &\leq & C\left( \varepsilon_2^{-1} + \int_{\varepsilon}^{\varepsilon_1} (q(\rho)  L(\rho) + \rho^{-1/2} L(\rho)^{1/2} + \rho q(\rho) +1) \, d\rho \right) \enspace .
\end{eqnarray*}
The first estimate follows from Lemma \ref{gronwall}. The second part is obtained, plugging the first estimate in the differential inequality (\ref{diffineq}). \ep

\begin{lem}\label{gronwall-2} Let $(s,k)\in [1,\infty)\times {\mathbb N}$ such that $s>k-1/2$. Suppose that the map defined on $(0,\varepsilon_0]$ by $\varepsilon \mapsto B(\varepsilon)$ is $C^1$ w.r.t the norm topology on ${\cal B}({\cal H})$ and that for any $\varepsilon \in (0,\varepsilon_0]$, $U$ and $B(\varepsilon)$ belong to $C^1(A)$. Assume that there exists $C>0$ such that for all $\varepsilon \in (0,\varepsilon_0]$, $q(\varepsilon) \leq C\varepsilon^{k-1}$. Then, for all $j\in \{1,\ldots,k\}$, there exist $C_j>0$ and $H_j\in L^1((0,\varepsilon_2])$ such that for all $(\varepsilon,z) \in (0,\varepsilon_2]\times \Omega_{(\theta_0 - \delta, \theta_0 + \delta),2}^+,$
\begin{eqnarray*}
\|F_{j,s}^{\pm}(\varepsilon,z)\| &<& C_j \\
\|\partial_{\varepsilon}F_{j,s}^{\pm}(\varepsilon,z)\| &\leq & H_j(\varepsilon) \enspace.
\end{eqnarray*}
\end{lem}
\noindent {\bf Proof:} Let us fix first $s\geq 1$. Since, it follows from the hypotheses that
\begin{equation*}
\int_0^{\varepsilon_0} q(\varepsilon)\, d\varepsilon < \infty \enspace,
\end{equation*}
the proof of the lemma for $j=1$ is a straightforward adaptation of the proof of Lemma \ref{gronwall-1}. In particular, $\sup_{(\varepsilon,z) \in (0,\varepsilon_2]\times \Omega_{(\theta_0 - \delta, \theta_0 + \delta),2}^+} \|F_{1,s}^{\pm}(\varepsilon,z) \| < \infty$. Applying Lemmata \label{lem4-1}, \ref{diffin2} and \ref{gronwall} implies the first estimate when $j\geq 2$. The second estimate is obtained, plugging the first estimate into the differential inequality (\ref{diffineq-2}). \ep

Let us explicit the implications of Lemmata \ref{gronwall-1} and \ref{gronwall-2}:
\begin{cor}\label{cor1} Under the hypotheses of Lemma \ref{gronwall-1}, we have that:
\begin{gather*}
\sup_{z \in S^+_{(\theta_0 - \delta, \theta_0 + \delta),\infty}} \|\bra A \ket^{-1}(1-zU^*)^{-1}\bra A \ket^{-1} \| < \infty\\
\sup_{z \in S^+_{(\theta_0 - \delta, \theta_0 + \delta),\infty}} \|\bra A \ket^{-1}(1-\bar{z}^{-1}U^*)^{-1}\bra A \ket^{-1} \| < \infty \enspace .
\end{gather*}
Let $\theta \in (\theta_0 - \delta, \theta_0 + \delta)$. Then,
\begin{itemize}
\item If $z$ tends to $e^{i\theta}$, then $\bra A \ket^{-1}(1-zU^*)^{-1}\bra A \ket^{-1}$ (resp. $\bra A \ket^{-1}(1-\bar{z}^{-1}U^*)^{-1}\bra A \ket^{-1}$) converges in norm (uniformly in $\theta$) to a bounded operator denoted $F_{1,1}^+(0^+,e^{i\theta})$ (resp. $F_{1,1}^-(0^+,e^{i\theta})$).
\item The maps defined on $(\theta_0 - \delta, \theta_0 + \delta)$ by $\theta \mapsto F_{1,1}^{\pm}(0^+,e^{i\theta})$ are norm-continuous functions with values in ${\cal B}({\cal H})$.
\end{itemize}
\end{cor}
\noindent{\bf Proof:} If $z\in {\mathbb D}^*$, the operators $T_{\varepsilon}^{\pm}(z)$ converge in norm respectively to $(1-zU^*)$ and $(1-\bar{z}^{-1}U^*)$ as $\varepsilon$ tends to $0$. This implies that $G_{\varepsilon}^{\pm}(z)$ converge also in norm respectively to $(1-zU^*)^{-1}$ and $(1-\bar{z}^{-1}U^*)^{-1}$ as $\varepsilon$ tends to $0$. Due to Lemma \ref{gronwall-1},
\begin{gather*}
\sup_{z \in S_{(\theta_0 - \delta, \theta_0 + \delta),2}^+} \|\bra A \ket^{-1}(1-zU^*)^{-1}\bra A \ket^{-1} \| < \infty\\
\sup_{z \in S_{(\theta_0 - \delta, \theta_0 + \delta),2}^+} \|\bra A \ket^{-1}(1-\bar{z}^{-1}U^*)^{-1}\bra A \ket^{-1} \| < \infty \enspace .
\end{gather*}
Since for $|z|\leq 1/2$, $\|\bra A \ket^{-1}(1-zU^*)^{-1}\bra A \ket^{-1} \|$ and $\|\bra A \ket^{-1}(1-\bar{z}^{-1}U^*)^{-1}\bra A \ket^{-1} \|$ are uniformly bounded, the first statement follows. The rest of the proof is similar in both cases, so we drop the superscript $\pm$ until the end. Due to Lemma \ref{gronwall-1}, for all $(\varepsilon,\mu,z) \in (0,\varepsilon_2]^2 \times \Omega_{(\theta_0 - \delta, \theta_0 + \delta),2}^+$, ($\varepsilon \leq \mu$),
\begin{equation}\label{integrability0}
\|F_{1,1}(\mu,z)-F_{1,1}(\varepsilon,z) \| \leq \int_{\varepsilon}^{\mu} \|\partial_{\rho} F_{1,1}(\rho,z)\|\, d\rho \leq \int_{\varepsilon}^{\mu} H(\rho)\, d\rho \enspace,
\end{equation}
where $H\in L^1(0,\varepsilon_2)$. This implies that $F_{1,1}(\varepsilon,z)$ converges in norm to a bounded operator when $\varepsilon$ tends to $0$ (uniformly in $z$, $z\in \Omega^+_{(\theta_0 - \delta, \theta_0 + \delta),2}$). The limit is denoted by $F_{1,1}(0^+,z)$. Of course, if $z\in S_{(\theta_0 - \delta, \theta_0 + \delta),2}^+$, $F_{1,1}^+(0^ +,z)=\bra A \ket^{-1}(1-zU^*)^{-1}\bra A \ket^{-1}$, $F_{1,1}^-(0^ +,z)=\bra A \ket^{-1}(1-\bar{z}^{-1}U^*)^{-1}\bra A \ket^{-1}$. For all $0\leq \varepsilon \leq \mu \leq \varepsilon_2$ and all $(z,z_0)\in (\Omega_{(\theta_0 - \delta, \theta_0 + \delta),2}^+)^2$, we have that:
\begin{equation*}\label{continuity}
\|F_{1,1}(\varepsilon,z)-F_{1,1}(0^+,z_0) \| \leq \| F_{1,1}(\varepsilon,z)-F_{1,1}(\mu,z) \|+\| F_{1,1}(\mu,z)-F_{1,1}(\mu,z_0) \|+\| F_{1,1}(\mu,z_0)-F_{1,1}(0^+,z_0) \| \enspace.
\end{equation*}
Using inequality (\ref{integrability0}), it follows that given $\delta'>0$, there exists $\varepsilon_3 \in (0,\varepsilon_2]$ such that for all $(\varepsilon,\mu)\in (0,\varepsilon_3]^2$ and all $(z,z_0)\in (\Omega_{(\theta_0 - \delta, \theta_0 + \delta),2}^+)^2,$
\begin{eqnarray*}
\| F_{1,1}(\mu,z_0)-F_{1,1}(0^+,z_0) \| &\leq & \delta'\quad \text{and}\\
\| F_{1,1}(\varepsilon,z)-F_{1,1}(\mu,z) \| &\leq & \delta' \enspace.
\end{eqnarray*}
Fix $\mu=\varepsilon_3$. The map $z\mapsto F_{1,1}(\varepsilon_3,z)$ is clearly norm-continuous on $(\Omega_{(\theta_0 - \delta, \theta_0 + \delta),2}^+)^2$ and there exists $\delta'' >0$, such that for all $z\in \Omega_{(\theta_0 - \delta, \theta_0 + \delta),2}^+$ with $|z-z_0|< \delta''$:
\begin{equation*}
\| F_{1,1}(\varepsilon_2,z)-F_{1,1}(\varepsilon_2,z_0) \| \leq \delta' \enspace.
\end{equation*}
Summing up, we have just proven that given $\delta' >0$, there exist $\varepsilon_3 \in (0,\varepsilon_2]$ and $\delta''>0$ such that for all $(\varepsilon,z) \in \Omega_{(\theta_0 - \delta, \theta_0 + \delta),2}^+ \times (0,\varepsilon_3]$ with $|z-z_0|< \delta''$
\begin{equation*}
\|F_{1,1}(\varepsilon,z)-F_{1,1}(0^+,z_0) \| \leq 3 \delta' \enspace .
\end{equation*}
In particular, if $\varepsilon$ vanishes, for all $z \in \Omega_{(\theta_0 - \delta, \theta_0 + \delta),2}^+$ such that $|z-z_0|< \delta''$, $\|F_{1,1}(0^+,z)-F_{1,1}(0^+,z_0) \| \leq \delta'$. If $z$ belongs to $S_{(\theta_0 - \delta, \theta_0 + \delta),2}^+$, this means that $\bra A \ket^{-1}(1-zU^*)^{-1}\bra A \ket^{-1}$ converges in norm to $F_{1,1}(0^+,z_0)$ when $z$ tends to $z_0$, $z_0 \in \Omega_{(\theta_0 - \delta, \theta_0 + \delta),2}^+$. If $z$ and $z_0$ belong to $\partial {\mathbb D}\cap \Omega_{(\theta_0 - \delta, \theta_0 + \delta),2}^+$, this means that the function defined on $(\theta_0 - \delta, \theta_0 + \delta)$ by $\theta \mapsto F_{1,1}(0^+,e^{i\theta})$ is continuous in norm. \ep

Now, we turn to the implications of Lemma \ref{gronwall-2}:
\begin{cor}\label{cor2} Under the hypotheses of Lemma \ref{gronwall-2}, we have that for all $j\in \{1,\ldots,k\},$
\begin{gather*}
\sup_{z \in S_{(\theta_0 - \delta, \theta_0 + \delta),\infty}^+} \|\bra A \ket^{-s}(1-zU^*)^{-j}\bra A \ket^{-s} \| < \infty\\
\sup_{z \in S_{(\theta_0 - \delta, \theta_0 + \delta),\infty}^+} \|\bra A \ket^{-s}(1-\bar{z}^{-1}U^*)^{-j}\bra A \ket^{-s} \| < \infty \enspace .
\end{gather*}
Let $\theta \in (\theta_0 - \delta, \theta_0 + \delta)$. Then,
\begin{itemize}
\item If $z$ tends to $e^{i\theta}$, then $\bra A \ket^{-s}(1-zU^*)^{-j}\bra A \ket^{-s}$ (resp. $\bra A \ket^{-s}(1-\bar{z}^{-1}U^*)^{-j}\bra A \ket^{-s}$) converges in norm to a bounded operator denoted $F^+_{j,s}(0^+,e^{i\theta})$ (resp. $F^-_{j,s}(0^+,e^{i\theta})$).
\item The maps defined on $(\theta_0 - \delta, \theta_0 + \delta)$ by $\theta \mapsto F_{j,s}^{\pm}(0^+,e^{i\theta})$ are norm-continuous functions with values in ${\cal B}({\cal H})$.
\item The maps defined on $(\theta_0 - \delta, \theta_0 + \delta)$ by $\theta \mapsto F^{\pm}_{1,s}$ are of class $C^{k-1}$ on $(\theta_0 - \delta, \theta_0 + \delta)$ with respect to the norm topology on ${\cal B}({\cal H})$.
\end{itemize}
\end{cor}
\noindent{\bf Proof:} The first part is an adaptation of the proof of Corollary \ref{cor1}, where  $F^{\pm}_{j,s}(\varepsilon,z)$ replaces $F^{\pm}_{1,1}(\varepsilon,z)$ and using Lemma \ref{gronwall-2} instead of Lemma \ref{gronwall-1}. We focus our attention on the last statement. As in the proof of Corollary \ref{cor1}, we drop the superscript $\pm$. For all $j\in \{1,\ldots,k\}$, all $(\varepsilon,\mu,z) \in (0,\varepsilon_2]^2 \times \Omega_{(\theta_0 - \delta, \theta_0 + \delta),2}^+$, ($\varepsilon \leq \mu$),
\begin{equation}\label{integrability}
\|F_{j,s}(\mu,z)-F_{j,s}(\varepsilon,z) \| \leq \int_{\varepsilon}^{\mu} \|\partial_{\rho} F_{j,s}(\rho,z)\|\, d\rho \leq \int_{\varepsilon}^{\mu} H_j(\rho)\, d\rho \enspace,
\end{equation}
where $H_j\in L^1(0,\varepsilon_2)$. This means that $F_{j,s}(\varepsilon,z)$ converges in norm to a bounded operator denoted by $F_{j,s}(0^+,z)$ when $\varepsilon$ tends to $0$, (uniformly in $z$, $z\in \Omega_{(\theta_0 - \delta, \theta_0 + \delta),2}^+$). It follows that $\bra A \ket^{-s}(1-zU^*)^{-j}\bra A \ket^{-s}$ (resp. $\bra A \ket^{-s}(1-\bar{z}^{-1}U^*)^{-j}\bra A \ket^{-s}$) converges in norm to $F_{j,s}(0^+,z_0)$ if $z$ tends to $z_0$, whenever $z_0 \in \Omega_{(\theta_0 - \delta, \theta_0 + \delta),2}^+$. Given $\varepsilon \in (0,\varepsilon_2]$, the operator-valued map $\theta \mapsto F_{j,s}(\varepsilon, e^{i\theta})$ is smooth on $(\theta_0 - \delta, \theta_0 + \delta)$ with respect to the norm topology on ${\cal B}({\cal H})$ and:
\begin{equation}\label{derivative}
\partial_{\theta} F_{j,s}(\varepsilon, e^{i\theta}) = i j \left(F_{j+1,s}(\varepsilon, e^{i\theta})-F_{j,s}(\varepsilon, e^{i\theta})\right) \enspace .
\end{equation}
If $j\in \{1,\ldots,k-1\}$, we already know that when $\varepsilon$ vanishes $F_{j+1,s}(\varepsilon, e^{i\theta})$ and $F_{j,s}(\varepsilon, e^{i\theta})$ converge in norm to $F_{j+1,s}(0^+, e^{i\theta})$ and $F_{j,s}(0^+, e^{i\theta})$ respectively (uniformly in $\theta$, $\theta \in (\theta_0 - \delta, \theta_0 + \delta)$). In particular, due to identity (\ref{derivative}), $\partial_{\theta} F_{j,s}(\varepsilon, e^{i\theta})$ converge in norm to $i j (F_{j+1,s}(0^+, e^{i\theta})-F_{j,s}(0^+, e^{i\theta}))$ (uniformly in $\theta$, $\theta \in (\theta_0 - \delta, \theta_0 + \delta)$). This means that the map $\theta \mapsto F_{j,s}(0^+, e^{i\theta})$ is differentiable on $(\theta_0 - \delta, \theta_0 + \delta)$, and
\begin{equation*}
\partial_{\theta} F_{j,s}(0^+, e^{i\theta}) = i j \left(F_{j+1,s}(0^+, e^{i\theta})-F_{j,s}(0^+, e^{i\theta})\right) \enspace .
\end{equation*}
If we take into account the first part of the corollary, the map $\theta \mapsto F_{j,s}(0^+, e^{i\theta})$ is of class $C^1$ on $(\theta_0 - \delta, \theta_0 + \delta)$, with respect to the norm topology on ${\cal B}({\cal H})$. The last part follows by induction on $m$, $m\in \{1,\ldots,k-1\}$ since for all $(\varepsilon,\theta) \in (0,\varepsilon_2] \times (\theta_0 - \delta, \theta_0 + \delta),$
\begin{equation*}
\partial^m_{\theta} F_{1,s}(\varepsilon, e^{i\theta}) = \sum_{j=1}^{m+1} a_j(m) F_{j,s}(\varepsilon, e^{i\theta}) \enspace ,
\end{equation*}
where the coefficients $(a_j(m))_{j\in \{1,\ldots,m+1\}}$ can be computed inductively: $a_1 (m+1)=-ia_1(m)$,
\begin{eqnarray*}
a_j(m+1) &=& i(j-1)a_{j-1}(m)-ija_j(m)\\
a_{m+2}(m+1) &=& i(m+1) a_{m+1}(m)\enspace.
\end{eqnarray*}
The map $\theta \mapsto F_{1,s}(0^+, e^{i\theta})$ is of class $C^{k-1}$ on $(\theta_0 - \delta, \theta_0 + \delta)$, with respect to the norm topology on ${\cal B}({\cal H})$ and for all $m\in \{1,\ldots,k-1\},$
\begin{equation*}
\partial^m_{\theta} F_{1,s}(0^+, e^{i\theta}) = \sum_{j=1}^{m+1} a_j(m) F_{j,s}(0^+, e^{i\theta}) \enspace .
\end{equation*} \ep

\noindent{\bf Remark:} It is also possible to prove that the maps $\theta \mapsto \partial^{k-1}_{\theta} F^{\pm}_{1,s}(0^+, e^{i\theta})$ are (norm-) H\"older continuous. This will be explicited in an upcoming work. We refer to \cite{jmp} Theorem 2.2 for the details. \\

Now, our main task consists in building a suitable family $(B(\varepsilon))$, which satisfies the hypotheses of Lemmata \ref{gronwall-1} and \ref{gronwall-2}. This is the purpose of the next paragraph.

\subsection{Properties of the function $q$}

The conclusions of Theorem \ref{thm1} will be drawn once established the relationships between the regularity properties of $U$, the family $(B(\varepsilon))$ and the properties of the function $q$.

\begin{lem}\label{qepsilon} Assume that $U$ is propagating w.r.t $A$ and there exists a family of uniformly bounded operators $(B(\varepsilon))_{\varepsilon \in (0,\varepsilon_0]}$ on ${\cal B}({\cal H})$ such that:
\begin{itemize}
\item $\lim_{\varepsilon \rightarrow 0}\|B(\varepsilon)-B_1\| =0$,
\item the map $\varepsilon \mapsto B(\varepsilon)$ is $C^1$ w.r.t the norm topology on ${\cal B}({\cal H})$
\item for any $\varepsilon \in (0,\varepsilon_0]$, $U$ and $B(\varepsilon)$ belong to $C^1(A)$,
\item the map $\varepsilon \mapsto \|\partial_{\varepsilon} B(\varepsilon)\|+\|\mathrm{ad}_A B(\varepsilon)\| + \varepsilon^{-1} \|B(\varepsilon)-B_1\|$ belongs to $L^1(0,\varepsilon_0)$.
\end{itemize}
Then, the function $q$ defined by: $q(\varepsilon) = \varepsilon^{-1} \max (\sup_{z \in \Omega_{{\mathbb T},2}^+}\|Q^{\pm}(\varepsilon,z)\|)$ for any $\varepsilon \in (0,\varepsilon_0]$ (with $Q^{\pm}(\varepsilon,z)$ defined in Lemma \ref{lem7}) belongs to $L^1(0,\varepsilon_0)$.
\end{lem}
\noindent{\bf Proof:} The conclusion follows from the definition of $Q^{\pm}(\varepsilon,z)$, Lemma \ref{lem7}, once noted that there exists $C>0$ such that for all $\varepsilon \in (0,\varepsilon_0]$: $\|e^{-\varepsilon B(\varepsilon)} \|\leq C$, $\|e^{\varepsilon B(\varepsilon)^*} \|\leq C$,
\begin{eqnarray*}
\| \sum_{k=1}^{\infty} \frac{(-\varepsilon)^k}{k!} \partial_{\varepsilon} (B(\varepsilon))^k \| &\leq & \sum_{k=1}^{\infty} \frac{\varepsilon^k}{k!} \| \partial_{\varepsilon} (B(\varepsilon))^k \| \leq C \varepsilon \| \partial_{\varepsilon} B(\varepsilon) \| \\
\|\mathrm{ad}_A e^{-\varepsilon B(\varepsilon)}\| &\leq & C \varepsilon \|\mathrm{ad}_A B(\varepsilon)\| \\
\|\mathrm{ad}_A e^{\varepsilon B(\varepsilon)^*}\| &\leq & C \varepsilon \|\mathrm{ad}_A B(\varepsilon)\|\enspace .
\end{eqnarray*}
\ep

\noindent{\bf Remark:} The construction of a family $(B(\varepsilon))_{\varepsilon \in (0,\varepsilon_0]}$ on ${\cal B}({\cal H})$ which satisfies the hypotheses of Lemma \ref{qepsilon} can be performed as soon as:
\begin{itemize}
\item $U\in C^1(A)$
\item $\mathrm{ad}_A U$ (or equivalently $U^*AU-A$) belongs to ${\cal C}^{0,1}(A)$.
\end{itemize}
See \cite{abmg} for the details.

Actually, the regularity properties of $U$ described in Theorem \ref{thm1} make possible the construction of a family $(B(\varepsilon))$ for which the function $q$ satisfies the hypotheses of Lemma \ref{gronwall-2} and Corollary \ref{cor2}. Let us introduce some local notations.

Given $k\in {\mathbb N}$, we set $\mathcal{N}=\{1,2,\ldots, k\}$. If for $j\in {\mathbb N}$, $\vec{\alpha}$ stands for the j-uple $(\alpha_1,\ldots, \alpha_j)$ in ${\cal N}^j$, we set: $|\vec{\alpha}| = \alpha_1 + \cdots +\alpha_j $ and $\vec{\alpha}\,! = \alpha_1 ! \, \cdots \, \alpha_j \,! $. For any collection of bounded operators in $\mathcal H$, $C_{\alpha_1},\ldots, C_{\alpha_j}$, we denote:
\begin{equation*}
\text{ad}_{C_{\vec{\alpha}}} (B) = \text{ad}_{C_{\alpha_1}}\circ \cdots \circ \text{ad}_{C_{\alpha_j}}(B) \enspace .
\end{equation*}

Let $U$ be a unitary operator which belongs to $C^{k+1}(A)$. Let us define the following sequence of bounded operators $B_p$, $p \in {\mathcal N}$ by: $B_1 = A - UAU^*$ and if $p \in {\cal N},$
\begin{align}
B_{p+1} = \, & p! \,\sum_{k=1}^p \frac{(-1)^k}{k!} \, \sum_{ \vec{\alpha} \in {\mathcal N}^k , |\vec{\alpha}|=p} \frac{1}{\vec{\alpha} !} \, \mathrm{ad}_{B_{\alpha_1}} \circ \cdots \circ \mathrm{ad}_{B_{\alpha_{k-2}}}(\mathrm{ad}_A B_{\alpha_{k-1}} ) \nonumber \\
& -p! \, \sum_{k=1}^p \frac{(-1)^{k+1}}{(k+1)!}   \sum_{ \vec{\alpha} \in {\mathcal N}^k , |\vec{\alpha}|+j=p+1 }\,
\frac{1}{\vec{\alpha} !\,  (j-1)!}\,  \mathrm{ad}_{B_{\vec{\alpha}}}(B_j) \enspace , \label{Qp}
\end{align}
This construction is actually motivated by Lemma \ref{npower}. A straightforward induction, based on Lemma \ref{equiv2}, leads us to the following lemma:
\begin{lem}\label{npower0} Assume that the unitary operator $U$ belongs to $C^{k+1}(A)$ for some $k\in {\mathbb N}$. Then for all $p\in \{1,2,\ldots, k+1\}$, $B_p \in C^{k+1-p}(A) $.
\end{lem}

Now, let us consider the family of bounded operators $(B(\varepsilon))_{\varepsilon \in (0,1]}$ defined by:
\begin{equation*}
B(\varepsilon) = \sum_{p=1}^k \,\frac{\varepsilon^{p-1}}{p!}\, B_p \enspace .
\end{equation*}
It follows that:
\begin{itemize}
\item the family $(B({\varepsilon}))$ is uniformly bounded,
\item $\lim_{\varepsilon \to 0} \| B(\varepsilon) - B_1 \|=0$,
\item the map $\varepsilon \to B(\varepsilon) $ is $C^1$ on $(0,1]$ with respect to the norm topology on ${\mathcal B}({\mathcal H})$,
\item for any $\varepsilon \in (0,1]$ , $B(\varepsilon) \in C^1(A)$.
\end{itemize}

Let us define the families $(Q^{\pm}(\varepsilon,z))_{(\varepsilon,z) \in (0,1]\times \Omega_{{\mathbb T},2}^+}$ by:
\begin{eqnarray*}
Q^+(\varepsilon,z) &=& z U^* \left( \partial_{\varepsilon} e^{-\varepsilon B(\varepsilon)} + B_1 e^{-\varepsilon B(\varepsilon)} - \mathrm{ad}_A e^{-\varepsilon B(\varepsilon)} \right) \\
Q^-(\varepsilon,z) &=& \bar{z}^{-1} U^* \left( \partial_{\varepsilon} e^{\varepsilon B(\varepsilon)^*} -B_1 e^{\varepsilon B(\varepsilon)^*} + \mathrm{ad}_A e^{\varepsilon B(\varepsilon)^*} \right) \enspace,
\end{eqnarray*}
which actually correspond to the definition of the families $(Q^{\pm}(\varepsilon,z))$ of Lemma \ref{lem7}. Then, we obtain that:
\begin{lem}\label{npower} Under the hypotheses of Lemma \ref{npower0}, there exists $M >0$ such that for all $\varepsilon \in (0,1]$, $q(\varepsilon) \leq M \varepsilon^k$, where $q(\varepsilon) = \varepsilon^{-1} \max (\sup_{z \in \Omega_{{\mathbb T},2}^+}\|Q^{\pm}(\varepsilon,z)\|)$ for any $\varepsilon \in (0,1]$.
\end{lem}
\noindent{\bf Proof:} Let us focus our attention on $\sup_{z \in \Omega_{{\mathbb T},2}^+}\|Q^+(\varepsilon,z)\|$. The other case can be treated similarly. There exists $M >0$ such that for all $(\varepsilon,z) \in (0,1]\times \Omega_{{\mathbb T},2}^+,$
\begin{equation*}
\|Q^+(\varepsilon,z)\| \leq M\, \|(\partial_{\varepsilon} e^{-\varepsilon B(\varepsilon)} ) e^{\varepsilon B(\varepsilon)} + B_1+ e^{-\varepsilon B(\varepsilon)}A
 e^{\varepsilon B(\varepsilon)}-A \| \enspace .
\end{equation*}
Since for any bounded operator $C$, $\|C\|=\|C^*\|$, the proof can be reduced to the control of the norm involved on the RHS of the first inequality. Writing  $C=\varepsilon B(\varepsilon)$, it follows from Corollary \ref{8-2} that for all $\varepsilon \in (0,1],$
\begin{eqnarray*}
e^{-C}A\,e^{C}-A & = & \sum_{k=1}^{\infty }  \frac{(-1)^{k+1}}{k!} \left( \sum_{\vec{\alpha}\in \mathcal{N}^k }\frac{ \varepsilon^{|\vec{\alpha}|}}{\vec{\alpha} !} \,\mathrm{ad}_{B_{\alpha_1}} \circ \cdots \circ \mathrm{ad}_{B_{\alpha_{k-2}}}(\mathrm{ad}_A B_{\alpha_{k-1}} )\right) \\
                 & = & \sum_{p=1}^{\infty} \varepsilon^p \left( \sum_{k=1}^p \frac{(-1)^{k+1}}{k!} \sum_{\vec{\alpha} \in \mathcal{N}^k,\,
|\vec{\alpha} | = p}    \frac{1}{\vec{\alpha}\, !} \, \mathrm{ad}_{B_{\alpha_1}} \circ \cdots \circ \mathrm{ad}_{B_{\alpha_{k-2}}}(\mathrm{ad}_A B_{\alpha_{k-1}} )\right)\enspace ,
\end{eqnarray*}
and
\begin{eqnarray*}
(\partial_{\varepsilon} e^{-C} ) e^{C} + B_1 & = &
 \sum_{p=1}^{\infty} \frac{(-1)^p}{p!}\,  \text{ad}_C^{p-1}(\partial_{\varepsilon}C)+B_1 \\
& = &  \sum_{p=2}^{\infty} \frac{ (-1)^p }{p!} \, \text{ad}_C^{p-1}(\partial_{\varepsilon}C)-  \sum_{j=1}^{k-1} \frac{\varepsilon^j}{j!}\, B_{j+1}\enspace .
\end{eqnarray*}
Recall that  $C$ stands for $\varepsilon B(\varepsilon) = \sum_{p=1}^{k} \frac{\varepsilon^p }{p!}\, B_p$. So, we get
\begin{eqnarray*}
  \sum_{k\geq 2} \frac{(-1)^k}{k!}\,  \text{ad}_C^{k-1}(\partial_{\varepsilon}C ) & =  &
   \sum_{k\geq 1} \frac{(-1)^{k+1}}{(k+1)!} \sum_{\vec{\alpha} \in {\mathcal N}^k}  \frac{\varepsilon^{|\vec{\alpha}|} }{\vec{\alpha} ! } \, \text{ad}_{B_{\vec{\alpha}}} (\partial_{\varepsilon} C) \, \\
    & = &
     \sum_{p=1}^{\infty} \varepsilon^p \left(
      \sum_{k=1}^p \frac{(-1)^{k+1}}{(k+1)!}   \sum_{ \vec{\alpha} \in {\mathcal N}^k , |\vec{\alpha}|+j=p+1 }\,
\frac{1}{\vec{\alpha} !\,  (j-1)!}\,  \text{ad}_{B_{\vec{\alpha}}}(B_j) \right) \enspace .
\end{eqnarray*}
Thus, the term $(\partial_{\varepsilon} e^{-\varepsilon \, B(\varepsilon)} ) e^{\varepsilon \, B(\varepsilon)} + (A-UAU^*)+ e^{-\varepsilon \, B(\varepsilon)}A
 e^{\varepsilon \, B(\varepsilon)}-A $ may be written as a norm convergent series of the form
 $\sum_{p=0}^{\infty} \varepsilon^p T^p$, where for all $p \leq k$,
 \begin{eqnarray*}
 T_p & = & \sum_{k=1}^p \frac{(-1)^{k+1}}{k!} \sum_{\vec{\alpha} \in {\mathcal N}^k , |\vec{\alpha}|=p} \frac{1}{\vec{\alpha}!}\, \mathrm{ad}_{B_{\alpha_1}} \circ \cdots \circ \mathrm{ad}_{B_{\alpha_{k-2}}}(\mathrm{ad}_A B_{\alpha_{k-1}} ) \\
     &   & +
      \sum_{k=1}^p \frac{(-1)^{k+1}}{(k+1)!}   \sum_{ \vec{\alpha} \in {\mathcal N}^k , |\vec{\alpha}|+j=p+1 }\,
\frac{1}{\vec{\alpha} !\,  (j-1)!}\,  \text{ad}_{B_{\vec{\alpha}}}(B_j) + \frac{1}{p!}\, B_{p+1}\\
& = & 0
 \end{eqnarray*}
due to our construction of the coefficients $(B_p)$ (See relation (\ref{Qp})). This finishes the proof. \ep

\subsection{Proof of Theorem \ref{thm1}}

Let us start with a consequence of Lemma \ref{qepsilon}:
\begin{prop}\label{thm0} Let $\Theta$ be an open subset of ${\mathbb T}$. Assume $U$ is propagating with respect to $A$ on $e^{i\Theta}$ and that there exist a family of uniformly bounded operators $(B(\varepsilon))_{\varepsilon \in (0,\varepsilon_0]}$ which satisfies the hypotheses of Lemma \ref{qepsilon}. Then, 
\begin{itemize}
\item[(i)] $U$ has a finite number of eigenvalues in $e^{i\Theta}$. Each of these eigenvalues has a finite multiplicity. The spectrum of $U$ has no singular continuous component in $e^{i\Theta}$.
\item[(ii)] For any compact subset $K\subset \Theta\setminus \sigma_{pp}(U),$
\begin{equation*}
\sup_{ |z|\neq 1, \arg z \in K} \|\bra A \ket^{-1}(1-zU^*)^{-1}\bra A \ket^{-1} \| < \infty \enspace .
\end{equation*}
\item[(iii)] If $z$ tends to $e^{i\theta}$, then $\bra A \ket^{-1}(1-zU^*)^{-1}\bra A \ket^{-1}$ converges in norm to a bounded operator denoted $F^+_{1,1}(0^+,e^{i\theta})$ (resp. $F^-_{1,1}(0^+,e^{i\theta}) )$ if $|z|<1$ (resp. $|z|>1$). This convergence is uniform if $\theta$ belongs to any compact subset $K\subset \Theta\setminus \sigma_{pp}(U)$.
\item[(iv)] The operator-valued functions defined by $F^{\pm}_{1,1}$ are continuous on $\Theta\setminus \sigma_{pp}(U)$, with respect to the norm topology on ${\cal B}({\cal H})$.
\end{itemize}
If $U$ is strictly propagating with respect to $A$ on $e^{i\Theta}$, then Statement (i) can be replaced by: $U$ is purely absolutely continuous on $e^{i\Theta}$.
\end{prop}
\noindent{\bf Proof:} In view of Corollary \ref{virial-2}, we know that $U$ has at most a finite number of eigenvalues in $\Theta$. These eigenvalues have finite multiplicity. We also know that for any $\theta \in \Theta\setminus \sigma_{pp}(U)$, there exists $\delta_{\theta}>0$ such that $U$ is strictly propagating with respect to $A$ on $(\theta-2\delta_{\theta},\theta+2\delta_{\theta})$. Given any compact subset $K\subset \Theta\setminus \sigma_{pp}(U)$, the collection $((\theta-\delta_{\theta},\theta+\delta_{\theta}))_{\theta \in \Theta\setminus \sigma_{pp}(U)}$ induces an open covering of $K$, from which we can extract a finite open covering. Due to Lemma \ref{qepsilon}, Corollary \ref{cor1} applies on each of these intervals, which proves all the statements of Theorem \ref{thm0}. In particular, we have that for any compact subset $K\subset \Theta\setminus \sigma_{pp}(U),$
\begin{gather*}
\sup_{z \in S_{K,\infty}^+} \|\bra A \ket^{-1}(1-zU^*)^{-1}\bra A \ket^{-1} \| < \infty\\
\sup_{z \in S_{K,\infty}^+} \|\bra A \ket^{-1}(1-\bar{z}^{-1}U^*)^{-1}\bra A \ket^{-1} \| < \infty \enspace .
\end{gather*}
This means that the operator $\bra A \ket^{-1} E(K)$ is $U$-smooth, which in turns implies that Ran$E(K) \subset {\cal H}_{ac}(U)$ (see \cite{abcf} Theorem 2.2 and Remark 3). Since $K$ was arbitrarily chosen in $\Theta\setminus \sigma_{pp}(U)$, the spectrum of $U$ has no singular component in $\Theta$. If in addition, $U$ is strictly propagating with respect to $A$ on $\Theta$, it clearly follows from Theorem \ref{virial} (or Theorem \ref{thm2}) that $U$ has even no eigenvalues in $\Theta$.
\\

\noindent{\bf Proof of Theorem \ref{thm1}:} Consider the family of uniformly bounded operators $(B(\varepsilon))_{\varepsilon \in (0,1]}$ defined by (\ref{Qp}). It satisfies the hypotheses of Lemma \ref{qepsilon} and Proposition \ref{thm0} due to Lemmata \ref{npower0} and \ref{npower}. In view of Corollary \ref{virial-2}, we know that $U$ has at most a finite number of eigenvalues in $\Theta$. These eigenvalues have finite multiplicity. We also know that for any $\theta \in \Theta\setminus \sigma_{pp}(U)$, there exists $\delta_{\theta}>0$ such that $U$ is strictly propagating with respect to $A$ on $(\theta-2\delta_{\theta},\theta+2\delta_{\theta})$. Given any compact subset $K\subset \Theta\setminus \sigma_{pp}(U)$, the collection $((\theta-\delta_{\theta},\theta+\delta_{\theta}))_{\theta \in \Theta\setminus \sigma_{pp}(U)}$ induces an open covering of $K$, from which we can extract a finite open covering. As noted above, due to Lemmata \ref{npower0} and \ref{npower}, Corollary \ref{cor2} applies on each of these intervals, which proves all the statements of Theorem \ref{thm1}. Since we have that for any compact subset $K\subset \Theta\setminus \sigma_{pp}(U),$
\begin{gather*}
\sup_{z \in S_{K,\infty}^+} \|\bra A \ket^{-1}(1-zU^*)^{-1}\bra A \ket^{-1} \| < \infty\\
\sup_{z \in S_{K,\infty}^+} \|\bra A \ket^{-1}(1-\bar{z}^{-1}U^*)^{-1}\bra A \ket^{-1} \| < \infty \enspace .
\end{gather*}
we also conclude as in the proof of Proposition \ref{thm0} that the spectrum of $U$ has no singular component in $\Theta$. Now, let $\phi \in C_0^{\infty}(\Theta\setminus \sigma_{pp}(U))$. It follows from the lines above that,
\begin{equation*}
\bra A\ket^{-s} \partial_{\theta}E(\theta) \bra A\ket^{-s} = \frac{1}{2\pi} \left(F^+_{1,s}(0^+,e^{i\theta}) - F^-_{1,s}(0^+,e^{i\theta})\right) \enspace .
\end{equation*}
Therefore, for all $m\in {\mathbb Z},$
\begin{equation*}
\bra A\ket^{-s} U^m \Phi(U) \bra A\ket^{-s} = \int_{\mathbb T} e^{i m\theta} \Phi(e^{i\theta})\, \bra A\ket^{-s} \partial_{\theta} E(\theta) \bra A\ket^{-s}\, d\theta \enspace .
\end{equation*}
If $s>k+1/2$, $\theta \mapsto \bra A\ket^{-s} \partial_{\theta} E(\theta) \bra A\ket^{-s}$ is of class $C^k$ with respect to the norm topology on the support of $\Phi$. Integration by part yields the existence of $C>0$ such that for all $m\in {\mathbb Z},$
\begin{equation*}
\| \bra A\ket^{-s} U^m \Phi(U) \bra A\ket^{-s}\| \leq C \bra m\ket^{-k} \enspace .
\end{equation*} \ep

\section{Regularity classes for Bounded Operators}

\subsection{Basics}

This section gathers some elementary properties of the regularity classes $C^k(A)$ (sometimes denoted $C^k(A,{\cal H})$ or $C^k(A,{\cal H},{\cal H})$) introduced in Section 2 and applies them in two specific contexts. For more details see \cite{abmg} Chapter 5. From now, $A$ is a fixed self-adjoint operator, densely defined on a fixed Hilbert space ${\cal H}$, with domain ${\cal D}(A)$. \\

The regularity of a bounded operator defined on ${\cal H}$ w.r.t $A$ is associated to the algebra of derivation on ${\cal B}({\cal H})$ defined by the operation $\mathrm{ad}_A$. From a theoretical point of view, it is often more convenient to reformulate this concept of derivation in terms of the regularity of the strongly continuous function:
\begin{eqnarray*}
{\cal W}_B: {\mathbb R} & \rightarrow & {\cal B}({\cal H})\\
t & \mapsto & e^{iAt}Be^{-iAt} \enspace .
\end{eqnarray*}
Most of the properties derived below can be deduced easily once established the following equivalence:
\begin{prop}\label{equiv0} Let $k\in {\mathbb N}$. The following assertions are equivalent:
\begin{itemize}
\item $B\in C^k(A)$
\item The map ${\cal W}_B$ is $C^k$ with respect to the strong topology on ${\cal B}({\cal H})$.
\item The map ${\cal W}_B$ is $C^k$ with respect to the weak topology on ${\cal B}({\cal H})$.
\end{itemize}
Moreover, ${\cal W}_B^{(k)} (0) = i^k \mathrm{ad}_A^k B$.
\end{prop}
For a proof, see \cite{abmg} Lemma 6.2.9, Theorem 6.2.10 in association with Lemma 6.2.1 and Definition 6.2.2. For all nonnegative integral number $k$, $C^{k+1}(A) \subset C^k(A)$.

\begin{prop}\label{2} If $B\in C^1(A)$, then $B({\cal D}(A)) \subset {\cal D}(A)$.
\end{prop}

For any nonnegative integral number $k$, $C^k(A)$ is clearly a vector subspace of ${\cal B}({\cal H})$. These classes also share the following algebraic properties:
\begin{prop}\label{3} Let $k\in {\mathbb N}$ and $(B,C) \in C^k(A)\times C^k(A)$. then,
\begin{itemize}
\item $B^* \in C^k (A)$ and for all $j\in \{0,\ldots,k\}$, $\mathrm{ad}_A^j B^* = (-1)^j (\mathrm{ad}_A^j B)^*$
\item $BC \in C^k(A)$ and for all $j\in \{1,\ldots,k\},$
\begin{equation*}
\mathrm{ad}_A^j BC = \sum_{l_1+l_2=j} \frac{j!}{l_1! l_2!} \mathrm{ad}_A^{l_1} B \mathrm{ad}_A^{l_2} C \enspace.
\end{equation*}
In particular, $\mathrm{ad}_A BC = (\mathrm{ad}_A B) C + B (\mathrm{ad}_A C)$
\item for all $j\in \{0,\ldots, k\}$, $\mathrm{ad}_A^j B \in C^{k-j}(A)$.
\item If $B$ is invertible (i.e $B^{-1}\in {\cal B}({\cal H})$) and $B\in C^1(A)$, then $B^{-1} \in C^1(A)$: $\mathrm{ad}_A B^{-1} = - B^{-1} (\mathrm{ad}_A B) B^{-1}$.
\end{itemize}
\end{prop}
See \cite{abmg} Propositions 5.1.2, 5.1.5, 5.1.6, 5.1.7 for a proof. Combining the last sentences of Proposition \ref{3}, we deduce that if an invertible bounded operator $B$ belongs to $C^k(A)$, then its inverse $B^{-1}$ also belongs to $C^k(A)$.

For the relationships between these regularity classes with the self-adjoint functional calculus, we refer the reader to \cite{abmg} Theorem 6.2.5 and Corollary 6.2.6 or to \cite{gj} and the Helffer-Sjostrand formula.


\subsection{Application to GGT Matrices}

In this paragraph, we show how the concepts of paragraph 7.1 can be implemented to measure the regularity of the GGT matrices considered in Section 5.

We follow the notations of Section 5: $(\gamma_k)$ and $(\beta_k)$ will stand for two sequences of ${\mathbb D}^{\mathbb Z}$ such that: $0<\inf_{k\in {\mathbb Z}} |\gamma_k|\leq \sup_{k\in {\mathbb Z}} |\gamma_k| < 1$ and $0<\inf_{k\in {\mathbb Z}} |\beta_k| \leq \inf_{k\in {\mathbb Z}} |\beta_k| < 1$. We will relate the regularity properties of the matrices $D_1(\gamma)$ and $H(\gamma)$.

The reader will note that given such a sequence $(\gamma_k)\in {\mathbb D}^{\mathbb Z}$, then there exists $\Phi\in C_0^{\infty}({\mathbb R}, {\mathbb R})$ with compact support in $(-\infty,1)$ such that:
\begin{equation*}
D_2(\gamma) = \sqrt{1-D_1(\gamma)^*D_1(\gamma)} =\Phi(D_1(\gamma)^*D_1(\gamma)) \enspace .
\end{equation*}
It follows that:
\begin{lem}\label{ggt1-0} If $D_1(\gamma)-D_1(\beta)$ is compact, so are $D_2(\gamma)-D_2(\beta)$ and $H(\gamma)-H(\beta)$.
\end{lem}
\noindent{\bf Proof:} The first part is a consequence of Stone-Weierstrass Theorem. In addition, we have that:
\begin{align}\label{ha-hb}
H(\gamma)-H(\beta) &= T(D_2(\gamma)-D_2(\beta)) - T^*(D_1(\gamma)-D_1(\beta))T (1-D_2(\gamma)T)^{-1} D_1(\gamma)^* \nonumber \\
& + T^*D_1(\beta))T (1-D_2(\gamma)T)^{-1} (D_2(\gamma)-D_2(\beta))T  (1-D_2(\gamma)T)^{-1}D_1(\gamma)^* \nonumber \\
& - T^*D_1(\beta) T (1-D_2(\beta)T)^{-1} (D_1(\gamma) - D_1(\beta))^*
\end{align}
Each term on the RHS is the product of bounded operators with at least one compact operator, which ends the proof. \ep

In other words, the difference $H(\gamma)-H(\beta)$ is compact whenever:
\begin{equation*}
\lim_{|k|\rightarrow \infty} (\gamma_k -\beta_k) = 0
\end{equation*}
In particular, due to Weyl's Theorem, we have that $\sigma_{\text{ess}}(H(\gamma)) = \sigma_{\text{ess}}(H(\beta))$.

Now, let us consider the regularity properties:
\begin{lem}\label{ggt2} Let $n\in {\mathbb N}$. If $D_1(\gamma)$ belongs to $C^n(B_a)$, then $D_2(\gamma)$ and $H(\gamma)$ also belong to $C^n(B_a)$.
\end{lem}
\noindent{\bf Proof:} The fact that $D_2(\gamma)$ belongs to $C^n(B_a)$ is a consequence of \cite{abmg} Theorem 6.2.5 and Corollary 6.2.6. Using Fourier transform or a direct proof by induction, we can show that $T$ and $T^*$ also belong to $C^{\infty}(B_a)$. Therefore, using the algebraic properties shared by the classes $C^n(B_a)$ (see e.g. Section 7.1 and references therein), the last part of the proof follows from formula (\ref{Uzero}). \ep

\begin{lem}\label{ggt2-1} Let $D_1(\gamma)$ and $D_1(\beta)$ in $C^1(B_a)$. Assume that the operators $(D_1(\gamma)-D_1(\beta))$ and $\mathrm{ad}_{B_a }(D_1(\gamma)-D_1(\beta))$ are compact. Then, $D_2(\gamma)-D_2(\beta)$ belongs to $C^1(B_a)$ and $\mathrm{ad}_{B_a} (D_2(\gamma)-D_2(\beta))$ is also compact.
\end{lem}
\noindent{\bf Proof:} We already know by Lemma \ref{ggt2} that $D_2(\gamma)-D_2(\beta)$ belongs to $C^1(B_a)$. Using sesquilinear forms and the Helffer-Sj\"ostrand formula (see e.g. \cite{D2}), $\mathrm{ad}_{B_a}(D_1(\gamma)-D_1(\beta))$ can be rewritten as a norm convergent integral of compact operators, so it is compact. \ep

\begin{lem}\label{ggt3} Let $D_1(\gamma)$, $D_1(\beta)$ in $C^1(B_a)$ such that $D_1(\gamma)-D_1(\beta)$ and $\mathrm{ad}_{B_a} (D_1(\gamma)-D_1(\beta))$ are compact. Then, the unitary operators $H(\gamma)$ and $H(\beta)$ belong to $C^1(B_a)$ and $\mathrm{ad}_{B_a} (H(\gamma)- H(\beta))$ is compact.
\end{lem}
\noindent{\bf Proof:} The first part follows from Lemma \ref{ggt2}. Due to identity (\ref{ha-hb}) and the hypotheses, $H(\gamma)-H(\beta)$ can be rewritten as a finite sum of the form:
\begin{equation*}
H(\gamma)-H(\beta) = \sum_j c_j X_{j,1}\ldots X_{j,q_j}
\end{equation*}
where $(c_j) \subset {\mathbb R}$ and the operators $(X_{j,m})$ belong to $C^1(B_a)$. It follows that:
\begin{equation}\label{ha-hb2}
\mathrm{ad}_{B_a} (H(\gamma)-H(\beta)) = \sum_j c_j \left[ (\mathrm{ad}_{B_a} X_{j,1})\ldots X_{j,q_j} + \ldots + X_{j,1}\ldots (\mathrm{ad}_{B_a} X_{j,q_j}) \right]
\end{equation}
Due to Lemma \ref{ggt2-1}, each terms on the RHS of (\ref{ha-hb2}) is a product of bounded operators with at least one compact factor. This implies our second claim. \ep

It remains to reinterpret this operator theoretic approach in the context of the Verblunsky coefficients. This will be made explicit in the next paragraph.

\subsection{The Diagonal Case}

Otherwise noted, we follow the notations of Section 5. If $\gamma:=(\gamma_k)_{k\in {\mathbb Z}}$ is a bounded sequence in ${\mathbb C}^{\mathbb Z}$, we denote by $D_{\gamma}$ the bounded linear operator defined by its action on the canonical orthonormal basis of $l^2({\mathbb Z})$ by: $D_{\gamma} e_k=\gamma_k e_k$. The letter $x$ stands for the sequence $(k)_{k\in {\mathbb Z}}$.

The family of seminorms $(p_{n_1,n_2})$ and $(q_n)$ are defined (for non-negative integral numbers $n_1$, $n_2$ and $n$) on ${\mathbb C}^{\mathbb Z}$ by:
\begin{equation*}
p_{n_1,n_2}(u)= \sup_{k\in {\mathbb Z}} | k^{n_1} (\Delta^{n_2} u)_k |
\end{equation*}
where $(\Delta u)_k = u_k -u_{k-1}$ for all $k\in {\mathbb Z}$ and
\begin{equation}\label{seminorm}
q_n(u)= \sum_{m= 0}^n p_{m,m} \enspace .
\end{equation}

For any bounded sequence $\gamma$, $D_{\gamma} \in C^{\infty}(A)$ and $\mathrm{ad}_A D_{\gamma}=0$. We also define the automorphism of ${\cal B}(l^2({\mathbb Z}))$, $\omega$ by: $\omega(D):= T^* DT$. If $S$ stands for the translation in ${\mathbb C}^{\mathbb Z}$ (i.e. $(S\gamma)_k:= \gamma_{k+1}$ for any $k\in {\mathbb Z}$), then for any bounded sequence $\gamma$, $\omega(D_{\gamma})= D_{S\gamma}$. Let $\alpha$, $\beta$ and $\gamma$ are three bounded sequences, we define the following bounded operators on $l^2({\mathbb Z})$: for $n\in {\mathbb Z}\setminus \{0\},$
\begin{equation*}
J_n(D_{\alpha},D_{\beta}):= T^n D_{\alpha} +D_{\beta} T^{-n} \enspace .
\end{equation*}
For $n\neq 0$, $J_n(D_{\alpha},D_{\beta}):=J_{-n}(\omega^{-n}(D_{\beta}),\omega^{-n}(D_{\alpha}))=J_{-n}(D_{S^{-n}\beta},D_{S^{-n}\alpha})$. The following lemmata are obtained by direct computations:
\begin{lem}\label{criteria0} Let $\alpha$ and $\beta$ be two bounded sequences in ${\mathbb C}^{\mathbb Z}$. For all $m\in {\mathbb Z}\setminus \{0\}$, $J_m(D_{\alpha},D_{\beta}) \in C^1(A)$ and
\begin{equation*}
\mathrm{ad}_A J_m(D_{\alpha},D_{\beta}) = m J_m(D_{\alpha},-D_{\beta}) = J_m(D_{m \alpha},D_{-m\beta})\enspace.
\end{equation*}
\end{lem}

\begin{lem}\label{criteria1-0} Let $(m,n)\in {\mathbb Z}\times {\mathbb Z}\setminus \{0\}$, $z\in {\mathbb C}\setminus \{0\}$, $\alpha$ and $\beta$ be two bounded sequences in ${\mathbb C}^{\mathbb Z}$ such that $\sup_k |k(\alpha_{k+n}-\alpha-k)|<\infty$ and $\sup_k |k(\beta_{k+n}-\beta_k)|<\infty$. Then, $D_{\alpha}$ and $J_m(D_{\alpha}, D_{\beta})$ belong to $C^1(zT^n A+\bar{z}AT^{-n})$ and
\begin{itemize}
\item $\mathrm{ad}_{zT^n A+\bar{z}AT^{-n}} D_{\alpha} = J_n(D_{zx(\alpha-S^n\alpha)},D_{\bar{z}x(S^n\alpha-\alpha)})$,
\item if $n=-m,$ \begin{align*}
\mathrm{ad}_{T^n A+AT^{-n}} J_m(D_{\alpha},D_{\beta}) & = D_{zx(\alpha-S^n\alpha)+zm\alpha+\bar{z}x(S^n\beta-\beta)-m\beta} \\
& + J_{2m}(D_{\bar{z}(x-n)(\alpha-S^{-n}\alpha)+\bar{z}m\alpha},D_{z(S^{-n}\beta-\beta)(x-n)-zm\beta}) \enspace ,
\end{align*}
\item if $n=m,$ \begin{align*}
\mathrm{ad}_{T^n A+AT^{-n}} J_m(D_{\alpha},D_{\beta}) &= J_{2n}(D_{zx(\alpha-S^n\alpha)+zm\alpha}, D_{\bar{z}x(S^n\beta-\beta)-m\beta}) \\
&+ D_{\bar{z}(x-n)(\alpha-S^{-n}\alpha)+\bar{z}m\alpha+z(S^{-n}\beta-\beta)(x-n)-zm\beta} \enspace ,
\end{align*}
\item Otherwise, \begin{align*}
\mathrm{ad}_{T^n A+AT^{-n}} J_m(D_{\alpha},D_{\beta}) &= J_{n+m}(D_{zx(\alpha-S^n\alpha)+zm\alpha},D_{\bar{z}x(S^n\beta-\beta)-m\beta}) \\
&+ J_{m-n}(D_{\bar{z}(x-n)(\alpha-S^{-n}\alpha)+\bar{z}m\alpha}, D_{z(S^{-n}\beta-\beta)(x-n)-zm\beta}) \enspace .
\end{align*}
\end{itemize}
\end{lem}

Let us recall that the self-adjoint operator $B_a$ defined in Section 5 can be rewritten as follows:
\begin{equation*}
B_a= 2a(TA+AT^*)+a(T+T^*)-4A \enspace.
\end{equation*}

Let $(\gamma_k)$ be a bounded sequence in ${\mathbb C}^{\mathbb Z}$. Since the condition $\sup_k |k(\gamma_{k+1} -\gamma_k)| < \infty$, means that the sequence $x(S\gamma-\gamma)$ is bounded or equivalently that $p_{1,1}(\gamma)<\infty$, we have that:
\begin{cor}\label{criteria1} Let $m\in {\mathbb N}$. If $q_m(\gamma)<\infty$, then $D(\gamma)\in C^m(B_a)$.
\end{cor}
\noindent{\bf Proof:} By induction, applying Lemmata \ref{criteria0} and \ref{criteria1-0} with $n=1$. \ep
\\

\section{A Complement on weakly positive commutators}

As mentioned in Section 2, the conclusions of Theorem \ref{thm2} can be strengthened under stronger regularity assumptions. We made it explicit in Theorem \ref{weak}, which is actually the counterpart of \cite{bm} Theorem 2.1 in our unitary setting. Paragraphs 8.2 and 8.3 are devoted to its proof. Applications will be considered in a forthcoming work.

Until the end of this paragraph, $A$ denotes a fixed self-adjoint operator with dense domain ${\cal D}(A)$ on some fixed Hilbert space ${\cal H}$. We assume also that the unitary operator $U$ is weakly propagating with respect to $A$. In this context, $B$ stands for $B=A-UAU^*>0$.

\subsection{Hypotheses}

The framework of this discussion is presented with details in \cite{abmg} Paragraph 6.3. We recall the following notions. By a Friedrichs couple, we mean a couple of Hilbert spaces $({\cal H}_1,{\cal H}_2)$ such that ${\cal H}_1$ is continuously and densely embedded in ${\cal H}_2$. We also say that:
\begin{defin} Let $(W(x))_{x\in {\mathbb R}^n}$, $n\in {\mathbb N}$, a family of bounded operators defined on some Hilbert space ${\cal H}$. $(W(x))_{x\in {\mathbb R}^n}$ is a $C^0$-group if:
\begin{itemize}
\item $W(0)=I$ and $W(x+y)=W(x)+W(y)$ for all $(x,y)\in {\mathbb R}^n\times {\mathbb R}^n$,
\item the mapping $W: {\mathbb R}^n\rightarrow {\cal B}({\cal H})$ is strongly continuous.
\end{itemize}
If in addition for all $x\in  {\mathbb R}^n$, $W(x)$ is unitary, $(W(x))_{x\in {\mathbb R}^n}$ is a unitary $C^0$-group.
\end{defin}

One of the $C^0$-group considered hereafter will be defined by: $W(t)=W_{\cal H}(t):=e^{itA}$ where $t\in {\mathbb R}$. We also introduce a couple of additional Hilbert spaces:
\begin{defin} ${\cal S}$ denotes the completion of ${\cal H}$ for the norm $\|f\|_{\cal S}= \bra f, B f\ket^{1/2}$. ${\cal S}^*$ is the completion of $B{\cal H}$ for the norm $\|f\|_{{\cal S}^*}= \bra f, B^{-1} f\ket^{1/2}$
\end{defin}
By identifying ${\cal H}$ with its adjoint ${\cal H}^*$ (Riesz Lemma), $({\cal H},{\cal S})$ and $({\cal S}^*,{\cal H})$ are Friedrichs couples. ${\cal S}$ and ${\cal S}^*$ stay in duality with respect to the scalar product of ${\cal H}$: each of them will be identified to the other's adjoint. As a result, the operator $B$ extends to a unitary operator from ${\cal S}$ to ${\cal S}^*$.

Now, if we assume that all the operators $(W(t))_{t\in {\mathbb R}}$ leave ${\cal S}^*$ invariant, then, by restriction, this induces a $C_0$-group in ${\cal S}^*$ and by duality another in ${\cal S}$ writing: $W_{\cal S}(t)=(W(-t)|_{{\cal S}^*})^*$. Strictly speaking, the generator of $W_{\cal K}$ should be denoted by $A_{\cal K}$ for ${\cal K}={\cal S}^*, {\cal H}$ or ${\cal S}$ according to the context and its corresponding domains by ${\cal D}(A,{\cal K})$. For simplicity the generators will be denoted $A$, but we shall distinguish their domains. If $T\in {\cal B}({\cal S}, {\cal S}^*)$, then for any $t\in {\mathbb R}$, $W_{{\cal S}^*}(-t)TW_{\cal S}$ also belongs to $ {\cal B}({\cal S}, {\cal S}^*)$. This leads us to the following extension of the regularity concept described in Section 7:
\begin{defin} Let $T\in {\cal B}({\cal S}, {\cal S}^*)$. We say that $T\in C^1(A;{\cal S}, {\cal S}^*)$ if the sesquilinear form $Q$ defined on ${\cal D}(A,{\cal S})\times {\cal D}(A,{\cal S})$ (equipped with the induced topology of ${\cal S}\times {\cal S}$) by
\begin{equation*}
Q(\varphi,\psi):=\langle A\varphi,B\psi\rangle - \langle \varphi,BA\psi\rangle
\end{equation*}
extends continuously as a bounded form on ${\cal S}\times {\cal S}$. The operator assigned to the extension of $Q$ is denoted $[A,T]\in {\cal B}({\cal S},{\cal S}^*)$.
\end{defin}
\noindent{\bf Remark:} Let $T\in {\cal B}({\cal S}, {\cal S}^*)$. Then $T\in C^1(A;{\cal S}, {\cal S}^*)$ if and only if the map with value in ${\cal B}({\cal S}, {\cal S}^*)$ defined by ${\cal W}_T: t\mapsto W_{{\cal S}^*}(-t)TW_{\cal S}$ is strongly $C^1$. In this case, the strong derivative can be computed as.
\begin{equation*}
{\cal W}_T' (0)= i [A,T]\enspace .
\end{equation*}

Lastly, if ${\cal D}(A,{\cal S}^*)$ is equipped with the Hilbert structure associated to the graph norm: $\|f\|_{{\cal D}(A,{\cal S}^*)} =(\|f\|_{{\cal S}^*}^2+\|Af\|_{{\cal S}^*}^2)^{1/2}$, then $({\cal D}(A,{\cal S}^*),{\cal S}^*)$ is another Friedrichs couple. This allows us to introduce the interpolation space ${\cal K}=({\cal D}(A,{\cal S}^*),{\cal S}^*)_{1/2,1}$, which is actually densely embedded in ${\cal S}^*$. For more details, we refer to \cite{trieb}, \cite{abmg} and \cite{bm}. We can formulate now the main result of this section:
\begin{thm}\label{weak} Assume that $U$ is weakly propagating with respect to $A$ (in the sense of definition \ref{propagating}) and that $A-UAU^*$ belongs to $C^1(A,{\cal S}, {\cal S}^*)$. Then, there exists $C>0$ such that for all $z\in {\mathbb C}\setminus {\mathbb S}$ and all $\varphi\in {\cal K},$
\begin{equation*}
|\bra \varphi,(1-zU^*)^{-1}\varphi \ket | \leq C \|\varphi\|^2_{\cal K}
\end{equation*}
In particular, $U$ is purely absolutely continuous.
\end{thm}
The next paragraphs are devoted to the proof of Theorem \ref{weak}.

\subsection{Differential inequalities}

The core of the proof is an avatar of Mourre's differential inequality strategy. Let us introduce some local notations: if $r>1,$
\begin{eqnarray*}
S_r^+ &=& \{z\in {\mathbb C}; r^{-1}<|z|<1\}\\
S_r^- &=& \{z\in {\mathbb C}; 1<|z|<r\}\enspace .
\end{eqnarray*}
For $\varepsilon >0$ and $z\in {\mathbb C}\setminus \{0\}$, we define:
\begin{eqnarray*}
T_{\varepsilon}^{+}(z) &=& 1- z U^* e^{-\varepsilon B}\\
T_{\varepsilon}^{-}(z) &=& 1- \bar{z}^{-1} U^* e^{\varepsilon B}
\end{eqnarray*}
The following observation will be used without any further comment: for any $\varepsilon >0$, $\|e^{\pm\varepsilon B}-1\| \leq \varepsilon \|B\| e^{\varepsilon \|B\|}$.

Now, we deal with the invertibility of the family of bounded operators $(T_{\varepsilon}^{\pm}(z))$. Before, let us make a couple of remarks. First, if $A$ is a bounded invertible operator on ${\cal H}$, then, $(1-A)$ is invertible if and only if $(1-(A^{-1})^*)$ is invertible. In this case,
\begin{equation*}
\Re ((1+A)(1-A)^{-1}) = 2 \Re ((1-A)^{-1}) -1= (1-A)^{-1} - (1-(A^{-1})^*)^{-1}\enspace.
\end{equation*}
Next, the functions $h_1$ and $h_2$ defined on ${\mathbb R}$ by: $h_1(0)=h_2(0)=2$ and
\begin{eqnarray*}
h_1(x) &=& x^{-1} (1-e^{-2x})\\
h_2(x) &=& x^{-1} (e^{2x}-1)
\end{eqnarray*}
for $x\neq 0$, are homeomorphisms from ${\mathbb R}$ onto $(0,\infty)$, respectively monotone decreasing and monotone increasing. Since for any $\varepsilon \in [0,1]$, $1-e^{-2 \varepsilon B_1}=\varepsilon B_1 h_1(\varepsilon B_1)$ and $e^{2 \varepsilon B_1}-1=\varepsilon B_1 h_2(\varepsilon B_1)$, we have that: $c_1 \varepsilon B_1 \leq 1-e^{-2 \varepsilon B_1}$ and $c_2 \varepsilon B_1 \leq e^{2 \varepsilon B_1} -1$ for some positive constant $c_1$ and $c_2$.

Having this in mind, we state the following lemma:
\begin{lem}\label{lem4} The linear operators $T_{\varepsilon}^{\pm}(z)$ are invertible in ${\cal B}({\cal H})$, provided $(\varepsilon,z) \in [0,1]\times S_2^+$. Denote by $G_{\varepsilon}^{\pm}(z)$ the respective inverse of $T_{\varepsilon}^{\pm}(z)$. Then, there exists $C>0$ such that:
\begin{itemize}
\item For all $(\varepsilon,z) \in [0,1]\times S_2^+$: $\|G_{\varepsilon}^{\pm}(z)\| \leq C (1-|z|^2)^{-1}$.
\item For all $(\varepsilon,z) \in (0,1]\times S_2^+$ and all $\psi \in {\cal H}\supset {\cal S}^*,$
\begin{equation}
\|G_{\varepsilon}^{\pm}(z) \psi\|_{\cal S} \leq  C \sqrt{\frac{|\bra \psi, \Re(G_{\varepsilon}^{\pm}(z)) \psi \ket |}{\varepsilon}}\enspace.
\end{equation}
and subsequently $\|G_{\varepsilon}^{\pm}(z)\|_{{\cal S}^* \rightarrow {\cal S}} \leq C \varepsilon^{-1}$.
\end{itemize}
\end{lem}
\noindent {\bf Proof:} Let $c=\min(c_1,c_2)$. We have that for all $(\varepsilon,z) \in [0,1]\times S_2^+$:
\begin{equation*}
c |z|^{\pm 2} \varepsilon \bra \psi, B\psi\ket \pm (1-|z|^{\pm 2}) \|\psi \|^2 \leq \pm \bra \psi, (|z|^{\pm 2} (1- e^{\mp 2 \varepsilon B}) + (1-|z|^{\pm 2})) \psi \ket
\end{equation*}
which readily implies that:
\begin{eqnarray*}
\frac{c}{2} |z|^2 \varepsilon \|\psi\|_{\cal S}^2+ (1-|z|^2) \|\psi \|^2 &\leq & \bra T_{\varepsilon}(z)^{+} \psi, \psi \ket + \bra \psi, \bar{z}e^{- \varepsilon B} U T_{\varepsilon}(z)^{+} \psi \ket \\
\frac{c}{2} |z|^{-2} \varepsilon \|\psi\|_{\cal S}^2 - (1-|z|^{-2}) \|\psi \|^2 &\leq & \bra T_{\varepsilon}(z)^{-} \psi, \psi \ket + \bra \psi, z^{-1} e^{\varepsilon B} U T_{\varepsilon}(z)^{-} \psi \ket \enspace.
\end{eqnarray*}
This shows that the operators  $T_{\varepsilon}^{\pm}(z)$ are injective. On the other hand, Ran $T_{\varepsilon}^{+}(z)=$ Ker $(T_{\varepsilon}^{+}(z)^*)^{\perp} ={\cal H}$ since $T_{\varepsilon}^{+}(z)^*=-\bar{z}T_{\varepsilon}^{-}(z)e^{-\varepsilon B(\varepsilon)} U$. Similarly Ran $T_{\varepsilon}^{-}(z)=$ Ker $(T_{\varepsilon}^{-}(z)^*)^{\perp} ={\cal H}$. This proves the first part of the lemma. Let $(\varepsilon,z) \in [0,1]\times S_2^+$, $\varphi \in {\cal H}$ and $\psi = G_{\varepsilon}^{+}(z) \varphi$. Setting $a=c/2$, we have that:
\begin{eqnarray*}
a |z|^{\pm 2} \varepsilon \|G_{\varepsilon}^{\pm}(z) \varphi\|_{\cal S}^2 \pm (1-|z|^{\pm 2}) \|G_{\varepsilon}^{\pm }(z) \varphi \|^2 &\leq & \bra \varphi, G_{\varepsilon}^{\pm}(z) \varphi \ket - \bra \varphi, G_{\varepsilon}^{\mp}(z) \varphi \ket \\
&\leq & 2 \Re (\bra \varphi, G_{\varepsilon}^{\pm}(z) \varphi \ket) \enspace,
\end{eqnarray*}
which implies in particular that for $(\varepsilon,z) \in [0,1]\times S_2^+$: $\pm (1-|z|^{\pm 2}) \|G_{\varepsilon}^{\pm} (z) \varphi \|^2 \leq 2 \|\varphi\| \|G_{\varepsilon}^{\pm} (z) \varphi \|$. The first estimate follows. We also deduce that for all $(\varepsilon,z) \in (0,1]\times S_2^+$:
\begin{eqnarray*}
a c |z|^{\pm 2} \varepsilon\, \|G_{\varepsilon}^{\pm}(z) \varphi \|_{\cal S}^2 \leq 2|\bra \varphi, \Re (G_{\varepsilon}^{\pm}(z)) \varphi \ket| \leq 2 \|\varphi \|_{{\cal S}^*} \|G_{\varepsilon}^{\pm}(z) \varphi \|_{\cal S}\enspace,
\end{eqnarray*}
which implies the last estimates. \ep

Note that the maps $\varepsilon \mapsto e^{\pm \varepsilon B}$ are $C^1$ on ${\mathbb R}$ w.r.t the norm topology of ${\cal B}({\cal H})$. It follows that:
\begin{lem}\label{lem5} For any fixed $z \in S_2^+$, the maps $\varepsilon \mapsto T_{\varepsilon}^{\pm}(z)$ and $\varepsilon \mapsto G_{\varepsilon}^{\pm}(z)$ are $C^1$ on $(0,1]$ with respect to the norm topology on ${\cal B}({\cal H})$ and for all $(\varepsilon,z)\in (0,1]\times S_2^+$:
\begin{eqnarray*}
\partial_{\varepsilon} G_{\varepsilon}^{+}(z) &=& -z G_{\varepsilon}^{+}(z) U^* B e^{-\varepsilon B} G_{\varepsilon}^{+}(z)\\
\partial_{\varepsilon} G_{\varepsilon}^{-}(z) &=& \bar{z}^{-1} G_{\varepsilon}^{-}(z) U^* B e^{\varepsilon B} G_{\varepsilon}^{-}(z) \enspace .
\end{eqnarray*}
\end{lem}
\noindent{\bf Proof:} The regularity of the maps $\varepsilon \mapsto T_{\varepsilon}^{\pm}(z)$ follows from the previous remark: in particular,
\begin{eqnarray*}
\partial_{\varepsilon} T_{\varepsilon}^{+}(z) &=& z U^* B e^{-\varepsilon B}\\
\partial_{\varepsilon} T_{\varepsilon}^{-}(z) &=& - \bar{z}^{-1} U^* B e^{\varepsilon B} \enspace .
\end{eqnarray*}
Now, dropping the superscript $\pm$, we observe that given $z\in S_2^+$, for all $(\rho,\varepsilon) \in (0,1]^2,$
\begin{equation*}
G_{\rho}(z)-G_{\varepsilon}(z) = G_{\rho}(z) \left( T_{\varepsilon}(z)-T_{\rho}(z) \right) G_{\varepsilon}(z)
\end{equation*}
Due to Lemma \ref{lem5}, $\|G_{\varepsilon}(z)\| \leq C (1-|z|^2)^{-1}$ for all $\varepsilon \in [0,1]$. Since the map $\varepsilon \mapsto T_{\varepsilon}(z)$ is $C^1$ on $(0,1]$ w.r.t the norm topology of ${\cal B}({\cal H})$, it follows that the map $\varepsilon \mapsto G_{\varepsilon}(z)$ is also $C^1$ on $(0,1]$ w.r.t the same topology and that:
$\partial_{\varepsilon}G_{\varepsilon}(z) = - G_{\varepsilon}(z) ( \partial_{\varepsilon}T_{\varepsilon}(z) ) G_{\varepsilon}(z)$. The conclusion follows. \ep

Note that $C^1(A,{\cal S},{\cal S}^*)\subset C^1(A,{\cal H})$. It follows that for any $\varepsilon \in {\mathbb R}$, $e^{\varepsilon B}$ also belongs to $C^1(A,{\cal H})$ (see e.g. Lemma \ref{8-2}). In particular, $e^{\varepsilon B}{\cal D}(A,{\cal H})\subset {\cal D}(A,{\cal H})$ (see Proposition \ref{2}). This gives sense to the following lemma:
\begin{lem}\label{qepsilon} The sesquilinear forms $(Q_{\varepsilon})_{\varepsilon \in [-1,1]}$ defined a priori on ${\cal D}(A,{\cal H})^2$ by:
\begin{equation*}
Q_{\varepsilon}(\varphi,\psi) = \bra A\varphi, e^{\varepsilon B} \psi \ket - \bra e^{\varepsilon B}\varphi,A\psi \ket
\end{equation*}
extend continuously as bounded forms on ${\cal H}\times {\cal H}$. In particular, for all $(\varphi, \psi)\in {\cal H}\times {\cal H}$ and all $\varepsilon \in [-1,1]$, $| Q_{\varepsilon}(\varphi,\psi) | \leq C\,|\varepsilon| \|[A,B]\|\|e^{\varepsilon B}\varphi\|\|\psi\|$ for some $C>0$. If $B\in C^1(A,{\cal S},{\cal S}^*)$ then, for all $(\varphi, \psi)\in {\cal H}\times {\cal H}$ and all $\varepsilon \in [-1,1],$
\begin{equation*}
| Q_{\varepsilon}(\varphi,\psi) | \leq C\, |\varepsilon| \|[A,B]\|_{{\cal S}\rightarrow {\cal S}^*}\|e^{\varepsilon B}\varphi\|_{\cal S}\|\psi\|_{\cal S} \enspace ,
\end{equation*}
for some $C>0$.
\end{lem}
\noindent{\bf Proof:} The first assertion expresses the fact that the operator $e^{\varepsilon B}$ also belongs to $C^1(A)$ for any $\varepsilon \in {\mathbb R}$. As mentioned in our preliminary remark, for all $(\varphi, \psi)\in {\cal D}(A,{\cal H})^2,$
\begin{equation*}
Q_{\varepsilon}(\varphi,\psi) = \int_0^{\varepsilon}\bra A e^{(\varepsilon -\mu) B}\varphi, B e^{\mu B} \psi \ket - \bra B e^{(\varepsilon -\mu) B}\varphi, A e^{\mu B} \psi \ket \, d\mu \enspace .
\end{equation*}
The conclusion follows from the hypotheses, once noted that the family $(e^{\varepsilon B})_{\varepsilon \in [-1,1]}$ is uniformly bounded in ${\cal B}({\cal H})$ and extends as a uniformly bounded family of linear operators of ${\cal B}({\cal S})$. \ep

Since $(e^{\varepsilon B})_{\varepsilon \in [-1,1]} \subset C^1(A,{\cal H})$ and $U\in C^1(A,{\cal H})$, the operators $T_{\varepsilon}^{\pm}(z)$ and $G_{\varepsilon}^{\pm}(z)$ belong to $C^1(A,{\cal H})$ for any $(\varepsilon,z) \in [0,1]\times S_2^+$ (see Proposition \ref{3}). With Lemma \ref{qepsilon} in mind, we obtain in particular that: for $(\varepsilon,z) \in (0,1]\times S_2^+$ and any $(\varphi, \psi)\in {\cal D}(A,{\cal S}^*)^2$ (${\cal D}(A,{\cal S}^*)\subset {\cal H}$),
\begin{eqnarray}
\bra \varphi, \partial_{\varepsilon}G_{\varepsilon}^+(z)\psi \ket &=& - \bra A\varphi, G_{\varepsilon}^+(z)\psi \ket + \bra (G_{\varepsilon}^+(z))^*\varphi, A\psi \ket + z Q_{-\varepsilon}(U(G_{\varepsilon}^+)^*\varphi,G_{\varepsilon}^+\psi) \nonumber \\
\bra \varphi, \partial_{\varepsilon}G_{\varepsilon}^-(z)\psi \ket &=& \bra A\varphi, G_{\varepsilon}^-(z)\psi \ket - \bra (G_{\varepsilon}^-(z))^*\varphi, A\psi \ket - \bar{z}^{-1} Q_{\varepsilon}(U(G_{\varepsilon}^-)^*\varphi,G_{\varepsilon}^-\psi) \label{prederivative}\enspace .
\end{eqnarray}

Let us introduce more notations. Given any family of vectors $(\varphi_{\varepsilon})_{\varepsilon \in (0,1]} \subset {\cal D}(A,{\cal S}^*)$ such that the map $\varepsilon \mapsto \varphi_{\varepsilon}$ is $C^1$ w.r.t the topology defined by the norm $\|\cdot \|_{{\cal S}^*}$ (and incidentally w.r.t the topology of ${\cal H}$), we define the complex-valued functions $F^{\pm}$ on $(0,1]\times S_2^+$ by:
\begin{equation*}
F^{\pm}(\varepsilon,z) = \bra \varphi_{\varepsilon}, G_{\varepsilon}^{\pm}(z) \varphi_{\varepsilon}\ket \enspace.
\end{equation*}
It follows that:
\begin{lem}\label{diffin} Suppose that $U$ and $B$ belong respectively to $C^1(A)$ and $C^1(A,{\cal S},{\cal S}^*)$. Then, for any $z\in S_2^+$, the maps $\varepsilon \mapsto F^{\pm}(\varepsilon,z)$ are of class $C^1$ on $(0,1]$ and:
\begin{eqnarray*}
\partial_{\varepsilon}F^+(\varepsilon,z) &=& \bra \partial_{\varepsilon}\varphi_{\varepsilon}-A\varphi_{\varepsilon}, G_{\varepsilon}^+(z)\varphi_{\varepsilon} \ket + \bra G_{\varepsilon}^+(z)^*\varphi_{\varepsilon}, \partial_{\varepsilon}\varphi_{\varepsilon}+A\varphi_{\varepsilon} \ket + z Q_{-\varepsilon}(U(G_{\varepsilon}^+(z))^*\varphi_{\varepsilon}, G_{\varepsilon}^+(z)\varphi_{\varepsilon}) \\
\partial_{\varepsilon}F^-(\varepsilon,z) &=& \bra \partial_{\varepsilon}\varphi_{\varepsilon}+A\varphi_{\varepsilon}, G_{\varepsilon}^-(z)\varphi_{\varepsilon} \ket + \bra G_{\varepsilon}^-(z)^*\varphi_{\varepsilon}, \partial_{\varepsilon}\varphi_{\varepsilon}-A\varphi_{\varepsilon} \ket - \bar{z}^{-1} Q_{\varepsilon}(U(G_{\varepsilon}^-(z))^*\varphi_{\varepsilon}, G_{\varepsilon}^-(z)\varphi_{\varepsilon}) 
\end{eqnarray*}
It follows that there exists $C>0$ such that for all $(\varepsilon,z) \in (0,1]\times S_2^+,$
\begin{equation}\label{diffineq}
|\partial_{\varepsilon} F^{\pm}(\varepsilon,z) | \leq C\, \left( \sqrt{|F^{\pm}(\varepsilon,z) F^{\mp}(\varepsilon,z) |}+l(\varepsilon) \|\varphi_{\varepsilon}\|_{{\cal S}^*}^2 \right)+ l(\varepsilon) \varepsilon^{-1/2}\left( \sqrt{|F^{\pm}(\varepsilon,z) |}+\sqrt{|F^{\mp}(\varepsilon,z) |} \right)
\end{equation}
where $l(\varepsilon)=\|\partial_{\varepsilon} \varphi_{\varepsilon}\|_{{\cal S}^*}+\|A \varphi_{\varepsilon}\|_{{\cal S}^*}$.
\end{lem}
\noindent {\bf Proof:} The first part follows from identities (\ref{prederivative}). Then, we can deduce inequalities (\ref{diffineq}), by using Lemmata \ref{lem4} and \ref{qepsilon}, noting that:
\begin{eqnarray*}
\bar{z}e^{-\varepsilon B}U (G_{\varepsilon}^+(z))^* &=&-G_{\varepsilon}^-(z)\\
z^{-1} e^{\varepsilon B}U (G_{\varepsilon}^-(z))^* &=&-G_{\varepsilon}^+(z) \enspace ,
\end{eqnarray*}
and $(G_{\varepsilon}^{\pm}(z))^*= 1- G_{\varepsilon}^{\mp}(z)$ (the injection of ${\cal S}^*$ in ${\cal S}$ is continuous). \ep

The next step consists in integrating the differential inequality of Lemma \ref{diffin}. In order to apply successfully Lemma \ref{gronwall}, let us choose $\varphi \in {\cal K}$ where the interpolation space ${\cal K}:=({\cal S}^*,{\cal D}(A,{\cal S}^*))_{1/2,1}$ is continuously and densely embedded in ${\cal S}^*$ (see \cite{trieb} or \cite{abmg} Chapter 2 for the notations). Note that ${\cal D}(A,{\cal S}^*)$ is endowed with the Hilbert space structure associated to the norm
\begin{equation*}
\|\cdot \|_{{\cal D}(A,{\cal S}^*)}= (\|\cdot \|^2_{{\cal S}^*}+ \|A\cdot\|^2_{{\cal S}^*})^{1/2}
\end{equation*}
For such a vector $\varphi$, there exists a family of vectors $(\varphi_{\varepsilon})_{\varepsilon \in (0,1]} \subset {\cal D}(A,{\cal S}^*)$ such that the map $\varepsilon \mapsto \varphi_{\varepsilon}$ is $C^1$ and $\lim_{\varepsilon \rightarrow 0^+}\varphi_{\varepsilon} =\varphi$, both w.r.t the topology on ${\cal S}^*$ (and incidentally w.r.t the topology of ${\cal H}$). Actually, this construction can be explicited:
\begin{equation}\label{phiepsilon}
\varphi_{\varepsilon}=\varepsilon^{-1}\int_0^{\varepsilon}W(\tau)\varphi\, d\tau \enspace,
\end{equation}
(see \cite{bm}). In this case, the function $\varepsilon \mapsto \varepsilon^{-1/2} l(\varepsilon)$ is integrable, which also implies the integrability of the function $l$. Since for all $\varepsilon \in (0, 1],$
\begin{equation*}
\|\varphi_{\varepsilon}\|_{{\cal S}^*} \leq \|\varphi_1\|_{{\cal S}^*}+ \int_{1}^{\varepsilon}\|\partial_{\tau}\varphi_{\tau}\|_{{\cal S}^*}\, d\tau \enspace ,
\end{equation*}
the functions $\varepsilon \mapsto \|\varphi_{\varepsilon}\|_{{\cal S}^*}$ and $\varepsilon \mapsto l(\varepsilon)\|\varphi_{\varepsilon}\|_{{\cal S}^*}$ are also integrable.
As a consequence, we obtain:
\begin{lem}\label{gronwall-1} Let $\varphi \in {\cal K}$ and fix a family of vectors $(\varphi_{\varepsilon})_{\varepsilon \in (0,1]} \subset {\cal D}(A, {\cal S}^*)$ such that the map $\varepsilon \mapsto \varphi_{\varepsilon}$ is $C^1$ and $\lim_{\varepsilon \rightarrow 0^+}\varphi_{\varepsilon} =\varphi$ (both w.r.t the topology on ${\cal S}^*$). Suppose that $U$ and $B$ belong to $C^1(A)$ and $C^1(A,{\cal S},{\cal S}^*)$ respectively. Then there exist $C>0$ and $H\in L^1((0,1])$ such that for all $(\varepsilon,z) \in (0,1]\times S_{2}^+,$
\begin{eqnarray*}
|F^{\pm}(\varepsilon,z)| &<& C \\
|\partial_{\varepsilon}F^{\pm}(\varepsilon,z)| &\leq & H(\varepsilon) \enspace.
\end{eqnarray*}
\end{lem}
\noindent {\bf Proof:} The reader will observe first that the function $\varepsilon \mapsto \varepsilon \|\varphi_{\varepsilon} \|_{{\cal S}^*}^2$ is integrable. Define, the auxiliary functions $K$ and $L$ by
\begin{eqnarray*}
K(\varepsilon,z) &=& | F^+(\varepsilon,z)| +| F^-(\varepsilon,z)| \\
L(\varepsilon) &=& \sup_{z\in S_{2}^+} K(\varepsilon,z)
\end{eqnarray*}
Up some adjustment of the constants, we have that for all $(\varepsilon,z) \in (0,1]\times S_{2}^+,$
\begin{eqnarray*}
\left| K(1,z) - K(\varepsilon,z) \right| &= & \left| | F^+(1,z)| - | F^+(\varepsilon,z)|+  | F^-(1,z)| -| F^-(\varepsilon,z)| \right| \\
&\leq & | F^+(1,z) - F^+(\varepsilon,z)|+  | F^-(1,z) - F^-(\varepsilon,z)| \\
& \leq & \int_{\varepsilon}^{1} | \partial_{\rho} F^+(\rho,z) |+| \partial_{\rho} F^-(\rho,z) | \, d\rho \\
& \leq & C \int_{\varepsilon}^{1} ( K(\rho,z) +l(\rho) \rho^{-1/2} K(\rho,z)^{1/2} +l(\rho) \|\varphi_{\rho}\|_{{\cal S}^*}) \, d\rho 
\end{eqnarray*}
using Lemma \ref{diffin} and the fact that: $| F^{\pm}(\varepsilon,z)| \leq K(\varepsilon,z)$. It follows from Lemma \ref{lem4} that for all $(\varepsilon,z) \in (0,1]\times S_{2}^+,$
\begin{eqnarray*}
K(\varepsilon,z) &\leq & K(1,z) + C \int_{\varepsilon}^{1} (K(\rho,z) +l(\rho)\rho^{-1/2} K(\rho,z)^{1/2} +l(\rho) \|\varphi_{\rho}\|_{{\cal S}^*}) \, d\rho \\
L(\varepsilon) &\leq & L(1) + C \int_{\varepsilon}^{1} (q(\rho)  L(\rho) +l(\rho)\rho^{-1/2} L(\rho)^{1/2} +l(\rho) \|\varphi_{\rho}\|_{{\cal S}^*}) \, d\rho \enspace .
\end{eqnarray*}
The first estimate follows from Lemma \ref{gronwall}. The second part is obtained, plugging the first estimate in the differential inequality (\ref{diffineq}). \ep

\subsection{Proof of Theorem \ref{weak}}

Due to Lemma \ref{gronwall-1}, we conclude that the limits $F_0^{\pm}:=\lim_{\varepsilon \rightarrow 0^+} F_{\varepsilon}^{\pm}$ exist and satisfy:
\begin{equation*}
F^{\pm}(0,z) \leq F^{\pm}(1,z) +C \left(\int_0^1 \frac{d\rho}{\sqrt \rho} \|\varphi_{\rho}'\|_{{\cal S}^*}+\|A\varphi_{\rho}\|_{{\cal S}^*})^2 \right)^2 \enspace.
\end{equation*}
Following \cite{bm}, we know that:
\begin{equation*}
\int_0^1 \frac{d\rho}{\sqrt \rho} \|\varphi_{\rho}'\|_{{\cal S}^*}+\|A\varphi_{\rho}\|_{{\cal S}^*})^2 \leq \|\varphi\|_{\cal K} \enspace.
\end{equation*}
Using Lemma \ref{lem4} and identity $(\ref{phiepsilon})$, we have that:
\begin{equation*}
|F^{\pm}(1,z)| \leq \|G_1^{\pm}(z)\|_{{\cal S}^* \rightarrow {\cal S}} \|\varphi_1\|_{{\cal S}^*}^2 \leq C \left( \int_0^1 \|W(t)\varphi\|_{{\cal S}^*}^2 dt\right)^2 \leq C \|\varphi\|_{{\cal S}^*}^2 \leq C \|\varphi\|_{\cal K}^2 \enspace.
\end{equation*}
Those remarks imply that: $|F_0^{\pm}(z)|\leq C \|f\|_{\cal K}^2$ for some positive constant $C$. It remains to prove that given $z$, $|z|\neq 1$:
\begin{eqnarray*}
\lim_{\varepsilon \rightarrow 0} F^{+}(\varepsilon, z) &=& \bra \varphi, (1-zU^*)^{-1}\varphi\ket \\
\lim_{\varepsilon \rightarrow 0} F^{-}(\varepsilon, z) &=& \bra \varphi, (1-\bar{z}U^*)^{-1}\varphi\ket \enspace .
\end{eqnarray*}
Let us justify the first limit. Given $z\in S_2^+,$
\begin{eqnarray*}
| F^{+}(\varepsilon, z) - \bra \varphi, (1-zU^*)^{-1}\varphi\ket | &\leq & \|\varphi-\varphi_{\varepsilon}\|\| \left(\|G_{\varepsilon}^+(z)\| +\|G_0^+(z)\|\right) \|\varphi_{\varepsilon}\| \\
& +&\|G_{\varepsilon}^+(z)-G_0^+(z)\|\|\varphi_{\varepsilon}\|^2
\end{eqnarray*}
Due to lemma \ref{lem4}, $\|G_{\varepsilon}^+(z)\|\leq (1-|z|^2)^{-1}$, which combined with the second resolvent identity entails:
\begin{eqnarray*}
\|G_{\varepsilon}^+(z)-G_0^+(z)\| &\leq & C \|G_{\varepsilon}^+(z)\| \|1-e^{-\varepsilon B}\| \|G_0^+(z)\|\\
&\leq & C \varepsilon (1-|z|^2)^{-2} \enspace .
\end{eqnarray*}
On the other hand, since ${\cal S}^*$ is continuously embedded in ${\cal H}$, $\lim_{\varepsilon \rightarrow 0^+} \|\varphi_{\varepsilon}-\varphi\|=0$ and the family $(\varphi_{\varepsilon})_{\varepsilon \in (0,1]}$ is bounded w.r.t the Hilbert norm of ${\cal H}$. 
This allows to conclude the first part. The second affirmation is a consequence of \cite{abc} and the fact that the closure of ${\cal S}$ in the Banach space ${\cal K}^*$ is actually a closed subspace of ${\cal K}^*$ (see \cite{bm} p.4). \ep
\\

\noindent{\bf Acknowledgments:} The authors thank J. Asch for informative discussions.


\end{document}